\documentclass[11pt]{amsart}
\usepackage{amssymb,latexsym}
\usepackage{amsthm}
\usepackage{amsmath}
\usepackage[alphabetic]{amsrefs}
\usepackage{ae}
\usepackage{mathrsfs}
\newtheorem{theorem}{Theorem}
\newtheorem{prop}[theorem]{Proposition}
\newtheorem{cor}[theorem]{Corollary}
\newtheorem{lem}[theorem]{Lemma}

\theoremstyle{remark}
\newtheorem{rem}[theorem]{Remark}
\newtheorem{defin}[theorem]{Definition}

\numberwithin{equation}{section}

\newcommand{\twist}[2]{\,_{#1}#2}

\newcommand{\SP}{\operatorname{Sp}}
\newcommand{\GL}{\operatorname{GL}}

\newcommand{\PGL}{\operatorname{PGL}}
\newcommand{\SL}{\operatorname{SL}}
\newcommand{\Gal}{\operatorname{Gal}}

\newcommand{\Iso}{\operatorname{Iso}}
\newcommand{\Trd}{\operatorname{Trd}}
\newcommand{\Nrd}{\operatorname{Nrd}}
\newcommand{\Mat}{\operatorname{Mat}}
\newcommand{\Skew}{\operatorname{Skew}}

\newcommand{\Hom}{\operatorname{Hom}}

\newcommand{\tensor}{\otimes}

\newcommand{\XX}{\mathscr{X}}

\newcommand{\powfield}[1]{( \hspace{-1.6pt} ( {#1} ) \hspace{-1.6pt} )}

\newcommand{\f}{\mathfrak{f}}

\newcommand\R{\mathbf{R}}    %        real field 
\newcommand\C{\mathbf{C}}    %        complex field
\newcommand\Z{\mathbf{Z}}    %        integers
\newcommand\F{\mathbf{F}}    %        finite fields
\newcommand\Q{\mathbf{Q}}    %        rational field
\newcommand\G{\mathbf{G}}    %        mult or additive group
\newcommand\Proj{\mathbf{P}}      %        projective space
\newcommand\Lie{\operatorname{Lie}}
\newcommand\Ad{\operatorname{Ad}}
\newcommand\ad{\operatorname{ad}}

\newcommand\lie[1]{\mathfrak{#1}}
\newcommand\glie{\lie{g}}
\newcommand\mlie{\lie{m}}

\newcommand\plie{\lie{p}}

\newcommand{\hlie}{\lie{h}}
\newcommand{\ulie}{\lie{u}}
\newcommand{\vlie}{\lie{v}}
\newcommand{\wlie}{\lie{w}}

\newcommand{\clie}{\lie{c}}

\newcommand{\zlie}{\lie{z}}

\newcommand{\normal}{\lhd}

\newcommand{\bwedge}{\textstyle{\bigwedge}}
\newcommand{\congruent}{\equiv} 
\newcommand{\iso}{\simeq}
\renewcommand{\char}{\text{char}}

\newcommand{\sems}{\operatorname{ss}}

\newcommand{\NN}{\mathcal{N}}
\newcommand{\UU}{\mathcal{U}}
\newcommand{\VV}{\mathcal{V}}
\newcommand{\WW}{\mathcal{W}}

\newcommand{\LL}{\mathcal{L}}
\newcommand{\rank}{\operatorname{rank}}

\newcommand{\mm}{\mathfrak{m}}

\newcommand{\Int}{\operatorname{Int}}

\newcommand{\Stab}{\operatorname{Stab}}
\newcommand{\norm}[1]{\|#1\|}
\newcommand{\abs}[1]{|#1|}

\newcommand{\sep}{{\operatorname{sep}}}

\newcommand{\X}{\mathcal{X}}

\newcommand{\cc}{C_c}

\DeclareMathAlphabet{\matheucal}{U}{eus}{m}{n}
\newcommand{\orbit}{\matheucal{O}}

\newcommand{\tr}{\operatorname{tr}}

\begin{document}

\author{George J. McNinch}
\address{(personal address) 1715 East LaSalle Ave \\
         South Bend, Indiana 46617 \\
         USA}
\email{mcninchg@member.ams.org}
\date{17 January, 2004} 

\title{Nilpotent orbits over ground fields of good characteristic}

\begin{abstract}
  Let $X$ be an $F$-rational nilpotent element in the Lie algebra of a
  connected and reductive group $G$ defined over the ground field $F$.
  Suppose that the Lie algebra has a non-degenerate invariant bilinear
  form.  We show that the unipotent radical of the centralizer of $X$
  is $F$-split. This property has several consequences.  When $F$ is
  complete with respect to a discrete valuation with either finite or
  algebraically closed residue field, we deduce a uniform proof that
  $G(F)$ has finitely many nilpotent orbits in $\glie(F)$.  When the
  residue field is finite, we obtain a proof that nilpotent orbital
  integrals converge. Under some further (fairly mild) assumptions on
  $G$, we prove convergence for arbitrary orbital integrals on the Lie
  algebra and on the group.  The convergence of orbital integrals in
  the case where $F$ has characteristic 0 was obtained by Deligne and
  Ranga Rao (1972).
\end{abstract}

\maketitle
\tableofcontents

\section{Introduction}

Let $F$ be a field and $k$ an algebraically closed extension field.
Denote by $G$ a connected, reductive group defined over the ground field
$F$, and suppose the characteristic of $F$ to be good for $G$ (see \S
\ref{sec:general-reductive}).

The geometric nilpotent orbits, i.e. the nilpotent orbits of the
$k$-points of $G$ in $\glie = \glie(k)$, are described by the
Bala-Carter theorem; this result was proved for all good primes by
Pommerening.  Jens C. Jantzen has recently written a set of notes
\cite{jantzen:Nilpotent} on the geometric nilpotent orbits of a
reductive group; we refer to these notes -- and to their references --
for background on many of the results mentioned in this introduction.

The study of arithmetic nilpotent orbits, i.e. the nilpotent orbits of
the group of rational points $G(F)$ on the $F$-vector space
$\glie(F)$, is more complicated for a general field $F$; the description
of these orbits depends on Galois cohomology.  One of the goals of this paper is to
better understand the arithmetic nilpotent orbits.

In characteristic 0, or in large positive characteristic, the
Bala-Carter theorem may be proved by appealing to
$\lie{sl}_2$-triples.  To obtain a proof in any good characteristic,
other techniques are required. Pommerening's proof eventually shows
(after some case analysis) that one can associate to any nilpotent $X
\in \glie$ a collection of cocharacters of $G$ with favorable
properties; see \cite{jantzen:Nilpotent}.  Any cocharacter of $G$ determines a
parabolic subgroup, and it is a crucial result that each cocharacter
associated with $X$ determines the same parabolic subgroup $P$, which
therefore depends only on $X$.

On the other hand, Premet has recently given a more conceptual proof
of the Bala-Carter theorem. From the point of view of geometric
invariant theory, a vector $X \in \glie$ is nilpotent precisely when
its orbit closure contains 0; such vectors are said to be unstable.
According to the Hilbert-Mumford criterion for instability there is a
cocharacter $\phi$ of $G$ such that $X$ is unstable relative to the
$\G_m$-action on $\glie$ corresponding to $\phi$. A more precise form
of the Hilbert-Mumford criterion was established by Kempf and by
Rousseau; it yields cocharacters $\phi$ for which $X$ is in a suitable
sense optimally unstable relative to $\phi$.  Premet \cite{premet-kr} exploited these
optimal cocharacters, together with an idea of Spaltenstein, to prove
the Bala-Carter theorem in good characteristic.

Our first goal in this paper is to relate the associated cocharacters
found by Pommerening with the optimal cocharacters found by Premet;
this is done in Theorem \ref{theorem:assoc} after some
preliminaries in \S \ref{sec:general-instability}. We find that an
associated cocharacter for $X$ is optimal. An optimal
cocharacter $\phi$ need not be associated to $X$, but it almost is if
$X$ is a weight vector for the torus $\phi(\G_m)$. In particular,
the cocharacters associated with $X$ determine the same parabolic
subgroup $P$ as the optimal cocharacters for $X$; $P$ is 
called the instability parabolic (or instability flag) of $X$.

In a more general setting, Kempf exploited an important uniqueness
property of optimal cocharacters to prove that the instability
parabolic attached to an unstable $F$-rational vector is defined over
$F$, in case $F$ is perfect. In order to handle the case of an
imperfect field in the special case of the adjoint representation, we
invert this argument here.  Since a maximal torus of $G$ has at most
one cocharacter associated to $X$, it suffices to find a maximal $F$-torus
having a cocharacter associated to $X$; the
rationality of the cocharacter then follows from Galois descent. We
find such a torus under some assumptions on the separability of
orbits; the assumption usually holds for all nilpotent orbits in good
characteristic, at least when $G$ is semisimple.  The exception to
keep in mind is the group $\SL_n$ with $n$ divisible by the
characteristic. 

Since $X$ has an $F$-cocharacter associated to it, we deduce
more-or-less immediately that: (1) the instability parabolic $P$
attached to $X$ is defined over $F$ -- this had already been proved by
the author using other techniques; (2) the unipotent radical $R_u(C)$
of the centralizer $C$ of $X$ is defined over $F$ and is $F$-split;
and (3) $C$ has a Levi decomposition over $F$. See Theorem
\ref{theorem:split-u-radical} and Corollary
\ref{cor:centralizer-levi-decomp} for the latter two assertions. When
$F$ is perfect, e.g.  when $\char F = 0$, all unipotent groups over
$F$ are split.  See Remark \ref{rem:non-split-example} for an example
of a non-split unipotent group.

In \S \ref{sec:finiteness}, we study the Galois cohomological
consequences of the fact that $R_u(C)$ is $F$-split.  Suppose that $F$
is complete with respect to a non-trivial discrete valuation, and that
the residue field is finite or algebraically closed. If each adjoint
nilpotent orbit is separable, we prove that there are finitely many
arithmetic nilpotent orbits.  Our finiteness result improves one
obtained by Morris \cite{Morris-rational-conjugacy}. In \emph{loc.
  cit.}, the finiteness was obtained for various forms of classical
groups in good characteristic, and it was obtained for a general
reductive group under the assumption $p > 4h-4$ where $h$ is the
Coxeter number of $G$ (note that by now the use of the term ``very
good prime'' in \emph{loc.  cit.} \S3.13 is non-standard).

Suppose that $\glie$ has a non-degenerate invariant bilinear form.
This property guarantees that each geometric nilpotent orbit is
separable.  When the residue field of the complete field $F$ is
finite, it also guarantees that the centralizer of $X$ in $G(F)$ is a
unimodular locally compact group, so that the $G(F)$-orbit of $X$
carries an invariant measure.  When $\char F = 0$, a result of Deligne
and Rao asserts that this measure is finite for compact subsets of
$\glie(F)$. In \S \ref{sec:orbital-int}, we adapt the Deligne-Rao
argument to the case where $\char F >0$. We first treat the case where
$X$ is nilpotent; here we need no additional assumptions.  We obtain
convergence for orbital integrals of unipotent conjugacy classes in
$G(F)$ by invoking a result of Bardsley and Richardson which
guarantees the existence of a ``logarithm-like'' map $G \to \glie$;
these methods require some fairly mild assumptions on $G$ which are valid, for
example, when $G$ is a Levi subgroup of a semisimple group in
very good characteristic. Finally, we obtain the covergence for
general adjoint orbits and conjugacy classes under a somewhat
stronger additional assumption on the characteristic which guarantees
that the Jordan decomposition is defined over $F$.

Our proof that $R_u(C)$ is $F$-split answers a question put to the
author by D. Kazhdan; I thank him for his interest. I also thank S.
DeBacker, S. Evens, J. C. Jantzen, R. Kottwitz, J-P. Serre, and T.
Springer for useful conversations and comments regarding the
manuscript.

\section{Generalities concerning reductive groups}
\label{sec:general-reductive}

Recall that a homomorphism of algebraic groups $\varphi:A \to B$ is
said to be an \emph{isogeny} if it is surjective and has finite
kernel. The isogeny $\varphi$ will be said to be a \emph{separable
  isogeny} if $d\varphi:\Lie(A) \to \Lie(B)$ is an isomorphism. The
reader might keep the following example in mind as she reads the
material in this section: for any $n \ge 1$, the isogeny
$\varphi:\SL_{np/k} \to \PGL_{np/k}$ is not separable in
characteristic $p$. 

Throughout this section, $G$ is a connected and reductive group
defined over the infinite ground field $F$ of characteristic $p$.  The
field $k$ is an algebraically closed extension field of $F$.

\subsection{Good primes}

We first define the notions of \emph{good} and \emph{very good} primes
for  $G$. For a more thorough discussion of these
notions, the reader is referred to
\citelist{\cite{springer-steinberg} \cite{Hum95}
  \cite{jantzen:Nilpotent}}.  The reductive group $G$ is
assumed defined over $F$.

If $G$ is quasisimple with root system $R$, the characteristic $p$ of
$k$ is said to be bad for $R$ in the following circumstances: $p=2$ is
bad whenever $R \not = A_r$, $p=3$ is bad if $R = G_2,F_4,E_r$, and
$p=5$ is bad if $R=E_8$.  Otherwise, $p$ is good.  [Here is a more
intrinsic definition of good prime: $p$ is good just in case it
divides no coefficient of the highest root in $R$].

If $p$ is good, then $p$ is said to be very good provided that either
$R$ is not of type $A_r$, or that $R=A_r$ and $r \not
\congruent -1 \pmod p$.

If $G$ is reductive, the isogeny theorem \cite{springer-LAG}*{Theorem
9.6.5} yields a -- not necessarily separable -- central isogeny
$\prod_i G_i \times T \to G$ where the $G_i$ are quasisimple and $T$
is a torus. The $G_i$ are uniquely determined by $G$ up to isogeny,
and $p$ is good (respectively very good) for $G$ if it is good
(respectively very good) for each $G_i$.

The notions of good and very good primes are geometric
in the sense that they depend only on $G$ over an algebraically closed
field. Moreover, they depend only on the isogeny class of the derived
group $(G,G)$.

A crucial fact is the following:
\begin{lem}
  \label{lem:quasisimple-simple-adjoint}
  Let $G$ be a quasisimple group in very good characteristic.  Then the
  adjoint representation of $G$ on $\Lie(G)$ is irreducible and self-dual.
\end{lem}

\begin{proof}
  See \cite{Hum95}*{0.13}.
\end{proof}

We also note:
\begin{lem}
  \label{lem:good-for-subgroups}
  Let $M \le G$ be a reductive subgroup containing a maximal torus of $G$.
  \begin{enumerate}
  \item If $p$ is good for $G$, then $p$ is good for $M$.
  \item Suppose that $p > \rank_{\sems} G +1$, where $\rank_{\sems} G$
    is the semisimple rank of $G_{/k}$ (the \emph{geometric}
    semisimple rank). Then $p$ is \emph{very good} for $M$.
  \end{enumerate}
\end{lem}

\begin{proof}
  For (1), see for instance \cite{sommers-mcninch}*{Prop. 16}.  Now an
  inspection when $G$ is quasisimple yields (2).
\end{proof}

\subsection{Standard hypotheses}
We now recall the following \emph{standard
  hypotheses} SH for $G$; cf. \cite{jantzen:Nilpotent}*{\S 2.9}:

\begin{enumerate}
\item[SH1] The derived group of $G$ is simply connected.
\item[SH2] The characteristic of $k$ is good for $G$.
\item[SH3] There exists a $G$-invariant non-degenerate bilinear form
  $\kappa$ on $\glie$.
\end{enumerate}

\begin{defin}
  The reductive group $G$ will be said to be \emph{standard} if there
  is a separable isogeny between $G$ and a reductive group $H$ which
  satisfies the standard hypotheses SH. If $G$ is an $F$-group, then
  $G$ is $F$-standard if at least one such isogeny is defined over $F$.
\end{defin}

Observe that SH3 is preserved under separable isogeny; thus any
standard group has a non-degenerate invariant form $\kappa$ on its Lie
algebra.  Moreover, if the standard group $G$ is defined over a ground
field $F$, we may (and will) suppose that $\kappa$ is defined over
$F$.  Indeed, $\kappa$ amounts to an isomorphism between the
$G$-modules $\glie$ and $\glie^\vee$, so we need to find such an
isomorphism defined over $F$. Since $F$ is infinite,
$\Hom_{G/F}(\glie(F),\glie^\vee(F))$ is a dense subset of
$\Hom_G(\glie,\glie^\vee)$ and so the (non-empty, Zariski open) subset
$\text{Isom}_G(\glie,\glie^\vee)$ has an $F$-rational point.

\begin{prop}
  \label{prop:standard-groups}
  \begin{enumerate}
  \item[(a)] A semisimple group in very good characteristic is standard.
  \item[(b)] The group $\GL(V)$ is standard for any finite
    dimensional vector space $V$.
  \item[(c)] If $G$ and $H$ are standard, then $G \times H$ is standard.
  \item[(d)] The centralizer $M$ of a semisimple element of a standard
    group is standard. Especially, a Levi subgroup of a standard group
    is standard.
  \end{enumerate}
\end{prop}

\begin{proof}  
  When $G$ is simply connected, Lemma
  \ref{lem:quasisimple-simple-adjoint} implies that SH holds for $G$.
  If $\pi:G_{\operatorname{sc}} \to G$ is the simply connected
  covering isogeny, that same lemma implies that $d\pi$ is an
  isomorphism, hence that $\pi$ is a \emph{separable} isogeny; this
  proves (a) in general.
  
  For (b), the only thing that needs verifying is SH3; for this, it is
  well-known that the trace form on $\lie{gl}(V)$ is non-degenerate.
  Assertion (c) is straightforward. 
  
  For (d), note that the characteristic is good for $M$ by Lemma
  \ref{lem:good-for-subgroups}. Let $\varphi$ be a separable isogeny
  between $G$ and a group $\widehat G$ satisfying SH. The Levi
  subgroup $M$ of $G$ is the connected centralizer of a suitable
  semisimple element $s \in G$; let $\hat s \in \widehat G$ correspond via
  $\varphi$ to $s$ (thus either $\hat s = \varphi(s)$ or $\hat s \in
  \varphi^{-1}(s)$).  Then $\widehat M = C^o_{\widehat G}(\hat s)$ is a
  Levi subgroup of $\widehat G$, and $\varphi$ restricts to a
  separable isogeny between $M$ and $\widehat M$. The assertion (d)
  now follows from the following lemma.
\end{proof}

\begin{lem}
  \label{lem:semisimpl-centralizer}
  Let SH3 hold for $G$, let $x \in G$ be semisimple, and let $M =
  C_G^o(x)$ be the connected centralizer of $x$. Then $\kappa$
  restricts to a non-degenerate form on $\mlie = \Lie(M)$.
\end{lem}

\begin{proof}
  Write $\glie = \bigoplus_{\gamma \in k^\times} \glie_\gamma$ where
  for each $\gamma \in k^\times$, $\glie_\gamma$ is the
  $\gamma$-eigenspace of the (diagonalizable) map $\Ad(x)$.  Since
  $\kappa$ is invariant, the restriction $\kappa:\glie_\gamma \times
  \glie_{1/\gamma} \to k$ is a perfect pairing. The lemma now follows
  since $\mlie = \glie_1$.
\end{proof}

Let $X \in \glie$ and $g \in G$. When
$G$ is standard, the orbits of $X$ and $g$ are reasonably behaved:
\begin{prop}
  \label{prop:standard=>separable}
  Assume that $G$ is standard. The geometric orbits of $X \in \glie$
  and $g \in G$ are separable.  In particular, if $X \in \glie(F)$ and
  $g \in G(F)$, the centralizers $C_G(X)$ and $C_G(g)$ are defined
  over $F$.
\end{prop}

\begin{proof}
  Apply \cite{springer-steinberg}*{I.5.2 and I.5.6} for the first
  assertion. The fact that the centralizers are defined over $F$ then
  follows from \cite{springer-LAG}*{Prop. 12.1.2}.
\end{proof}
In general, of course, the $G$-orbit of an $F$-rational element which
is not separable need not be defined over $F$; such as orbit is
defined over $F$ if and only if it is defined over a separable closure
of $F$.

\section{The instability parabolic and nilpotent orbits}
\label{sec:general-instability}

In this section, we are concerned with a connected, reductive group
$G$ over an algebraically closed field $k$ whose
characteristic is good for $G$.

As described in the introduction, our goal here is to relate the
constructions given by Premet \cite{premet-kr} in his recent
simplification of the Bala-Carter-Pommerening Theorem to constructions
described in Jantzen's recent notes
\cite{jantzen:Nilpotent}. The main
result is Theorem \ref{theorem:assoc}.

\subsection{Length and cocharacters of $G$}
\label{sub:length}

Fix for a moment a maximal torus $T$ of $G$, and consider the lattice
$X_*(T)$ of cocharacters of $T$. Fix a $W$-invariant positive
definite, bilinear form $\beta$ on $X_*(T) \tensor \Q$. Given any
other torus $T' < G$, we may write $T' = \Int(g)T$, and one gets by
transport of structure a $W$-invariant form $\beta'$ on $X_*(T')$.
Since $\beta$ is $W$-invariant, $\beta'$ is independent of the choice
of the element $g$ with $T'=\Int(g)T$.

The form $\beta$ being fixed, there is a
unique $G$-invariant function $(\phi \mapsto \norm{\phi}):X_*(G) \to
\R_{\ge0}$ with the property $\norm \phi = \sqrt{\beta(\phi,\phi)}$
for $\phi \in X_*(T)$.

By a \emph{length function} $\norm\cdot$ on $X_*(G)$, we mean a
$G$-invariant function $\phi \mapsto \norm{\phi}$ associated with some
positive definite bilinear form $\beta$ on $X_*(T) \tensor \Q$ for
some maximal torus $T$ of $G$ in the above sense.  For the most part,
the choice of $T$ and $\beta$ will be fixed and we will not refer to
it.

%% Hmm. this seems to be unnecessary.
%%
%%
%Now suppose that $G$ is defined over the ground field $F$.
%Let $F_\sep$ denote the separable closure of $F$ in $k$, and
%let $\Gamma$ denote the Galois group $\gal(F_\sep/F)$. 
%Then $\Gamma$ acts on $X_*(G)$ as follows:
%for $\gamma \in \Gamma$, $\psi \in X_*(G)$ and $t \in F_\sep$,
%one has
%\begin{equation}
%  \label{eq:gal-on-cochars}
%  (\gamma \cdot \psi)(t) = \gamma(\psi(\gamma^{-1}(t))).
%\end{equation}
%We may suppose that $T$ is defined over $F$.  Then $\Gamma$ acts on
%$X_*(T)$. Denote by $\widehat W$ the subgroup of $\GL(X_*(T))$
%generated by the Weyl group $W = N_G(T)/T$ and the (finite) image of
%$\Gamma$.  Since the action of $\Gamma$ normalizes $W$, $\widehat W$
%is a finite group and so we may suppose that $\beta$ is $\widehat W$
%invariant (by ``averaging''). Thus, we may suppose, as in
%\cite{Kempf-instability}*{\S4}, that $\norm\psi = \norm{\gamma \cdot
%  \psi}$ for each $\psi \in X_*(G)$ and each $\gamma \in \Gamma$.

For later use, we observe the following:
\begin{lem}
  \label{lem:length-isogeny}
  Suppose $\pi:G \to G'$ is a surjective homomorphism of reductive
  groups with central kernel.  If $\norm{\cdot}'$ is a given length
  function on $X_*(G')$, there is a length function $\norm\cdot$ on
  $X_*(G)$ such that $\norm{\pi \circ \phi}' = \norm\phi$ for all
  $\phi \in X_*(G)$.
\end{lem}
\begin{proof}
  Fix a maximal torus $S'$ of $G'$ and a positive definite form
  $\beta'$ on $X_*(S') \tensor \Q$ giving rise to $\norm\cdot'$.
  Under our assumptions on $\pi$, $S=\pi^{-1}S'$ is a maximal torus of
  $G$. Since $\pi_{\mid S}:S \to S'$ is a surjective map of tori, the
  image of the induced map $X_*(S) \to X_*(S')$ has finite index in
  $X_*(S')$. Moreover, $\pi$ induces an isomorphism on Weyl groups
  $W=N_G(S)/S \xrightarrow{\sim} W'=N_{G'}(S')/S'$.  Now choose a
  $W$-module splitting of the exact sequence
  \begin{equation*}
    X_*(S) \tensor \Q \to X_*(S') \tensor \Q \to 0
  \end{equation*}
  so that $X_*(S) \tensor \Q \iso K \oplus (X_*(S')\tensor \Q)$ for
  some $W$-submodule $K$. Let $\beta''$ be a positive definite
  $W$-invariant bilinear form on $K$, and let $\beta = \beta'' \oplus
  \beta'$ be the corresponding form on $X_*(S) \tensor \Q$. One may
  then construct the length function on $X_*(G)$ using
  $\beta$, and the desired property is evident.
\end{proof}

A similar  observation is:
\begin{lem}
  \label{lem:length-subgroup}
  Let $G \subset G'$ be reductive groups and suppose that $G$ contains
  a maximal torus of $(G',G')$.  If $\norm\cdot$ is a given length
  function on $X_*(G)$, one can choose a length function $\norm\cdot'$
  on $X_*(G')$ such that $\norm\phi' = \norm\phi$ for $\phi \in
  X_*(G)$.
\end{lem}

\begin{proof}[Sketch]
  Let $S$ be a maximal torus of $G$ containing a maximal torus of
  $(G',G')$, and suppose $S \subset S'$, with $S'$ a maximal torus of
  $G'$.  Then $N_G(S)$ normalizes $S'$, and the map $N_G(S)/S \to
  N_{G'}(S')/S'$ is injective. The proof is now similar to that of the
  previous lemma.
\end{proof}

\subsection{Cocharacters and parabolic subgroups}
\label{sub:cochars-and-parabolics}

Let $\phi$ be a cocharacter of $G$.  Let 
\begin{equation*}
  P=P(\phi) = \{g \in G \mid \lim_{t\to 0} \phi(t)g\phi(t^{-1}) 
  \text{ exists}\}.
\end{equation*}
Then $P$ is a
parabolic subgroup of $G$
\cite{springer-LAG}*{8.4.5}. For $i \in \Z$, let
$\glie(i) = \glie(i;\phi)$ be the $i$-th weight space for
$\phi(\G_m)$. Then
\begin{equation*}
  \Lie(P) = \plie=
\bigoplus_{i \ge 0} \glie(i)
\end{equation*} 
The unipotent radical of $P$ is $U=\{g\in P \mid \lim_{t\to 0}
\phi(t)g\phi(t^{-1}) = 1\}$; see \cite{springer-LAG}*{8.4.6 exerc.
5}.  We have $\Lie(U) = \ulie =
\bigoplus_{i>0} \glie(i)$.

Note that $P(\phi) = P(n\phi)$ for any $n \in \Z_{\ge 1}$.

\begin{lem}
  \label{lem:phi-no-fixed-points}
  If $P=P(\phi)$ has unipotent radical $U$, then under the conjugation
  action, the torus $\phi(\G_m)$ has no fixed points $\not = 1$ on $U$.
\end{lem}
\begin{proof}
  This is immediate from the above description of $U$.
\end{proof}

\begin{lem}
  \label{lem:parabolic-homog}
  \begin{enumerate}
  \item Let $X \in \glie(i;\phi) \subset \ulie$ for $i \ge 1$,
    and let $u \in U$.  Then
    \begin{equation*}
      \Ad(u)X = X + \sum_{j>i} X_j \quad \text{for } X_j \in \glie(j;\phi).
    \end{equation*}
  \item Let $u \in U$ and put $\psi = \Int(u) \circ \phi$.  If $X \in
    \glie(i;\phi) \cap \glie(j;\psi)$ for some $i,j \ge 1$, then $i=j$ and $u \in
    C_P(X)$.
  \end{enumerate}
\end{lem}

\begin{proof}
  For (1), choose a maximal torus $T$ of $P$ containing $\phi(\G_m)$,
  and choose a Borel subgroup $B$ of $P$ containing $T$. Let $R
  \subset X^*(T)$ be the roots, let $R^+ \subset R$ be the positive
  system of non-zero $T$-weights on $\Lie(B)$, and let $R_U \subset R^+$
  be the $T$-weights on $\Lie(U)$; thus $R_U$ consists of
  those $\alpha \in R$ with $\langle \alpha,\phi \rangle > 0$.
  
  Let the homomorphisms $\XX_\alpha:\G_a \to B$ parameterize the root
  subgroups corresponding to $\alpha \in R$.  Then as a variety, $U$
  is the product of the images of the $\XX_\alpha$ for $\alpha \in
  R_U$. So to prove (1), it suffices to suppose $u = \XX_\beta(t)$ for
  $t \in \G_a$ and $\beta \in R_U$.
  
  Since $\glie(i;\phi) = \sum_{\alpha \in R_U; \langle \alpha,\phi
    \rangle = i} \glie_\alpha$ for $i \ge 1$, it suffices to suppose
  that $X \in \glie_\alpha$ with $\langle \alpha,\phi \rangle = i >
  0$.  In fact, we may suppose that $X = d\XX_\alpha(1)$.  By the
  \emph{Steinberg relations} \cite{springer-LAG}*{Prop 8.2.3} we have
  for $s \in \G_a$
  $$u\XX_\alpha(s)u^{-1} = \XX_\alpha(s) \cdot \prod_{\gamma \in R_U;
    \langle \gamma, \phi \rangle > i} \XX_{\gamma}(c_\gamma(s))$$
  for
  certain polynomial functions $c_\gamma(s)$.  Differentiating
  this formula, we get
  $$\Ad(u)X = X + \sum_{\gamma \in R_U; \langle \gamma, \phi \rangle >
    i} X_\gamma$$ with $X_\gamma \in \glie_\gamma$ as desired.
  
  For (2), note that since $X \in \glie(j;\psi)$ we have $\Ad(u^{-1})X
  \in \glie(j;\Int(u^{-1}) \circ\psi) = \glie(j;\phi).$ On the other
  hand, since  by assumption $X \in \glie(i;\phi)$, (1) shows
  that
  $$\Ad(u^{-1})X = X + \sum_{l > i} X_l$$
  with $X_l \in
  \glie(l;\phi)$.  Since the component of $\Int(u^{-1})X$ of weight
  $i$ for $\phi(\G_m)$ is the non-zero vector $X \in \glie(i;\phi)$,
  and since $\Int(u^{-1})X \in \glie(j;\phi)$, we see that $i=j$, and
  $\Int(u^{-1})X = X$.
\end{proof}

\begin{rem}
  \label{rem:rational-parabolic}
  Suppose the reductive group $G$ is defined over the ground field
  $F$. If $\phi:\G_m \to G$ is a cocharacter defined over $F$, then
  $P(\phi)$ is an $F$-parabolic subgroup. Conversely, if $P$ is an
  $F$-parabolic subgroup, then $P = P(\phi)$ for some cocharacter
  $\phi$ defined over $F$; for these assertions, see
 \cite{springer-LAG}*{Lemma  15.1.2}.
\end{rem}

\begin{rem}
  If the cocharacter $\phi:\G_m \to G$ is non-trivial, and
  $\phi(\G_m)$ is contained in the derived group of $G$, then
  $P(\phi)$ is a \emph{proper} parabolic subgroup (indeed, the assumption
  means that $\langle \alpha, \phi \rangle < 0$ for some root
  $\alpha$, hence $\glie_\alpha \not \subseteq \plie(\phi)$).
\end{rem}

\subsection{Geometric invariant theory and optimal cocharacters}

If $(\rho,V)$ is a rational representation of $G$, a vector $0\not=v \in V$
is unstable if the orbit closure $\overline{\rho(G)v}$ contains $0$.
For an unstable $v$ and a cocharacter $\phi$ of $G$, write $v =
\sum_{i \in \Z} v_i$, where $v_i \in V(i;\phi)$ and $V(i;\phi)$ is the $i$-th weight
space for $\phi(\G_m)$. We now define
\begin{equation*}
  \mu(v,\phi) = \min\{i \in \Z \mid v_i \not = 0\}.
\end{equation*}

A co-character $\phi$ of $G$ is optimal for an unstable vector $v \in
V$ if
\begin{equation*}
  \mu(v,\phi)/\norm\phi \ge \mu(v,\psi)/\norm \psi
\end{equation*}
for each cocharacter $\psi$ of $G$. This notion of course depends on the choice
of the length function $\norm{\cdot}$ on $X_*(G)$.
A co-character $\phi \in X_*(G)$ is primitive if there is no 
$(\psi,n) \in X_*(G) \times \Z_{\ge 2}$ with $\phi = n\psi$.

\begin{prop}(Kempf \cite{kempf-instab}, Rousseau)
  \label{prop:kr}
  Let $(\rho,V)$ be a rational representation of the reductive group
  $G$, and let $0 \not = v \in V$ be an unstable vector. 
  \begin{enumerate}
  \item The function $\phi \mapsto \mu(v,\phi)/\norm\phi$ on the set
    $X_*(G)$ attains a maximum value $B$; the cocharacters $\phi$ with
    $\mu(v,\phi)/\norm\phi = B$ are the optimal cocharacters for $v$.
  \item If $\phi$ and $\psi$ are optimal cocharacters for $v$, then
    $P(\phi)=P(\psi)$.
  \end{enumerate}
  Let $P$ be the common parabolic subgroup of (2). We then have:
  \begin{enumerate}
  \item[(3)] Let $\phi$ be an optimal cocharacter for $v$. For each $x
    \in P$, the cocharacter $\phi' = \Int(x) \circ \phi$ is optimal
    for $v$. Conversely, if $\phi$ and $\phi'$ are primitive optimal
    cocharacters for $v$, then $\phi$ and $\phi'$ are conjugate under
    $P$.
  \item[(4)] For each  maximal torus $T$ of $P$, there is a unique
    primitive $\phi \in X_*(T)$ which is optimal for $v$.
  \item[(5)] We have $\Stab_G(v) \subset P.$
  \end{enumerate}
\end{prop}

Write $P(v)$ for the parabolic subgroup of part (2) of the
Proposition; it is known as the \emph{instability flag} or the
\emph{instability parabolic.}

\subsection{Optimal cocharacters and central surjections}

Let $\pi:G \to G'$ be a surjective homomorphism between reductive
groups with central kernel, and construct the length functions
$\norm\cdot$ on $X_*(G)$ and $\norm\cdot'$ on $X_*(G')$ as in Lemma
\ref{lem:length-isogeny}.  Let $(\rho,V)$ and $(\rho',V')$ be rational
representations of $G$ and $G'$ respectively, and let $f:V \to V'$ be
a $G$-module homomorphism (for the pull-back $G$-module structure on
$V'$). Suppose that every non-0 vector of $\ker f$ is semistable (i.e.
not unstable). [We will consider precisely this setup in the proof of
Proposition \ref{prop:premet} below].
\begin{lem}
  \label{lem:isogeny-optimal}
  If $0 \not = v \in V$ is unstable, then $\phi \in X_*(G)$ is optimal
  for $v$ if and only if $\phi'=\pi \circ \phi \in X_*(G')$ is optimal
  for $f(v)$.
\end{lem}

\begin{proof}
  Let $B$ be the maximal value of $\mu(v,\psi)/\norm\psi$ for $\psi
  \in X_*(G)$, and let $B'$ be the maximal value of
  $\mu'(f(v),\nu)/\norm\nu'$ for $\nu \in X_*(G')$.
  
  With $\phi,\phi'$ as in the statement of the lemma, notice that we
  may find a $\phi(\G_m)$-submodule $W$ of $V$ such that $V \iso \ker
  f \oplus W$ as $\phi(\G_m)$-modules. In particular, $f$ induces an
  isomorphism $f_{\mid W}:W \to f(V)$. By hypothesis, $\phi(\G_m)$
  acts trivially on $\ker f$.  If $\phi$ is optimal for $v$, then $v =
  \sum_{i>0} v_i$ with each $v_i \in W(i;\phi)$.  Then $f(v) =
  \sum_{i>0} f(v_i)$ and it is clear that $\mu(v,\phi) =
  \mu'(f(v),\phi')$.  Conversely, if $\phi'$ is optimal for $f(v)$,
  write $f(v) = \sum_{i > 0} x_i$ with $x_i \in V'(i;\phi')$.  Then
  $v$ may be uniquely written $\sum_{i>0} y_i + z$ for certain $y_i
  \in W(i;\phi)$ with $f(y_i) = x_i$ and $z \in \ker f$.  Since $v$ is
  unstable, we must have $z=0$ and it is clear again that $\mu(v,\phi)
  = \mu'(f(v),\phi')$.
  
  Moreover, we have $\norm\phi = \norm{\phi'}'$. So the result will
  follow if we show that $B=B'$.  By what was said above, we know that
  $B \le B'$.  To show that equality holds, choose $\gamma \in
  X_*(G')$ which is optimal for $f(v)$. Since $\pi$ is surjective with
  central kernel, there is, as in the proof of Lemma
  \ref{lem:length-isogeny}, an $n \in \Z_{\ge 1}$ and $\phi \in
  X_*(G)$ with $n\gamma = \pi \circ \phi$.  Then applying the
  preceding considerations to $\phi'=n\gamma$ we get
  \begin{equation*}
    B'=\mu'(f(v),\gamma)/\norm\gamma' = \mu'(f(v),n\gamma)/\norm{n\gamma}'
    =\mu(v,\phi)/\norm{\phi} \le B.
  \end{equation*}
  Thus $B=B'$ and the lemma follows.
\end{proof}

\begin{rem}
  With notations as in the previous lemma, $\phi'$ may fail to be
  primitive when $\phi$ is primitive.
\end{rem}

\subsection{Optimal cocharacters for nilpotent elements.}

We are going to describe here a recent  result of Premet
giving a new approach to the classification of
nilpotent orbits for $G$ in good characteristic.

The first thing to notice is the following: for the adjoint
representation of $G$, the unstable vectors are precisely the
nilpotent elements. Indeed, that $0$ lies in the closure of each
nilpotent orbit is a consequence of the finiteness of the number of
nilpotent orbits; see \cite{jantzen:Nilpotent}*{\S2.10}.  On the other
hand, let $\chi:\glie \to \mathbf{A}^r$ be the \emph{adjoint quotient
  map}; cf. \cite{jantzen:Nilpotent}*{\S7.12, 7.13}.  The fiber
$\chi^{-1}(\mathbf{0})$ is precisely $\NN$; see \textit{loc. cit.}
Proposition 7.13. If $X \in \glie$ is not nilpotent, then it is
contained in a fiber $\chi^{-1}(b)$ with $b \not = \mathbf{0}$; since
this fiber is closed and $G$-invariant, $0 \not \in
\overline{\Ad(G)X}$. This proves our observation. Given $X \in \glie$
nilpotent, the result of Kempf and Rousseau (Proposition
\ref{prop:kr}) yields optimal cocharacters for $X$, and Premet
\cite{premet-kr} used this fact to give a simple proof of the
Bala-Carter-Pommerening Theorem.

To discuss Premet's work, we must recall some terminology.  A
nilpotent $X \in \glie$ is said to be distinguished provided that the
connected center of $G$ is a maximal torus of $C_G(X)$.  A parabolic
subgroup $P <G$ is distinguished if
\begin{equation*}
  \dim P/U = \dim U/(U,U) + \dim Z
\end{equation*}
 where $U$
is the unipotent radical of $P$, and $Z$ is the center of $G$.

For $X \in \glie$ nilpotent, write $C=C_G(X)$ for the 
centralizer of $X$, and $P=P(X) = P(\phi)$ for
the instability parabolic subgroup, where
$\phi$ is any optimal cocharacter for $X$. Moreover,
write $U$ for the unipotent radical of $P$.

\begin{prop} (Premet)
  \label{prop:premet}
  Fix a length function on $X_*(G)$.
  If $X \in \glie$ is nilpotent, there
  is a cocharacter $\phi$ which is optimal for $X$ with the following
  properties:
 \begin{enumerate}
 \item $X \in \glie(2;\phi)$.
 \item The centralizer $C_\phi$ of $\phi(\G_m)$ in $C$ is reductive,
   and $C = C_\phi \cdot R$ is a Levi decomposition, where $R = C \cap
   U = R_u(C)$.
 \item Choose a maximal torus $S$ of $C_\phi$, and let $L = C_G(S)$.
   Then $X$ is distinguished in $\Lie(L)$, $P_L=P(\phi) \cap L$ is a
   distinguished parabolic subgroup of $L$, $X$ is in the open (Richardson)
   orbit of $P_L$ on its unipotent radical $U_L = U \cap L$, and
   $\phi(\G_m)$ lies in the derived subgroup of $L$.
 \end{enumerate}
\end{prop}

Note that in general neither $C$ nor $C_\phi$ is connected; the assertion
in (2) that $C_\phi$ is reductive is equivalent to: $C_\phi^o$ is reductive.
This proposition was proved by Premet \cite{premet-kr}*{Theorem 2.3,
  Proposition 2.5, Theorem 2.7} under the additional assumption that
$G$ satisfies the \emph{standard hypotheses} SH1--3 of \S
\ref{sec:general-reductive}.  Premet used the validity of the result
for this more restrictive class of groups $G$ to deduce a proof of the
Bala-Carter-Pommerening Theorem for any reductive group $G$ in good
characteristic. We will check here that the proposition itself is always true
in good characteristic.

\begin{proof}[Proof of Proposition \ref{prop:premet}]  
  Write $\norm{\cdot}_G$ for the fixed length function on $X_*(G)$.
  
  By \cite{springer-LAG}*{9.6.5}, we may find a central isogeny $\pi:H \to G$
  where $H = T \times \prod_i G_i$ and each $G_i$ is a simply
  connected, quasisimple group in good characteristic.  Since the
  characteristic is good, it follows from \cite{Hum95}*{0.13} that
  each proper $H$ submodule of $\bigoplus_i \Lie(G_i)$ is central in
  $\Lie(H)$.  Thus $\ker d\pi$ is central, and so $d\pi$ induces a
  bijection $\NN_H \to \NN_G$ by \cite{jantzen:Nilpotent}*{\S2.7};
  here, $\NN_G$ denotes the nilpotent variety of $G$, and $\NN_H$ that
  of $H$. We get also that each non-0 vector in $\ker d\pi$ is a
  semisimple element of $\glie$, hence is semistable (in fact:
  stable).
  
  We may choose a length function $\norm\cdot_H$ on $X_*(H)$
  compatible with $\norm{\cdot}_G$ as in Lemma
  \ref{lem:length-isogeny}.  We claim now that if the proposition
  holds for $H$ with this choice of length function, then it holds for
  $G$. To prove this claim, let $X \in \NN_G$, let $X' \in
  d\pi^{-1}(X)$ be the unique nilpotent preimage of $X$ in $\Lie(H)$,
  and let $\phi' \in X_*(H)$ satisfy the conclusion of the proposition
  for $H$. Put $\phi = \pi \circ \phi' \in X_*(G).$ We will show that
  $\phi$ satisfies the conclusion of the proposition for $G$.
  Property (1) needs no comment.  For (2), the only thing that must be
  verified is that $C_\phi$ is reductive.  Since $\pi$ restricts to a
  central isogeny $C'=C_H(X') \to C_G(X)$, it also restricts to a
  central isogeny $C'_{\phi'} \to C_\phi$. Since $C'_{\phi'}$ is
  reductive, $C_\phi$ is reductive, and (2) follows. Since $\pi$
  restricts to a central isogeny $P(\phi') \to P(\phi)$, the proof of
  (3) is similar.  It only remains to see that $\phi'$ is optimal for
  $X$; in view of our choice of $\norm{\cdot}_H$, this follows from Lemma
  \ref{lem:isogeny-optimal}. Our claim is proved.
  
  Finally, we may find a reductive group $M$ satisfying SH1--3 and an
  inclusion $H \subset M$ with $(M,M) = (H,H)$.  Indeed, for each $i$
  such that $G_i = \SL_n$, let $G_i' = \GL_n$, otherwise let
  $G_i'=G_i$; then take $M = T \times \prod G_i'$ with the obvious
  inclusion $H \subset M$.  As has already been remarked, Premet
  proved the proposition for $M$ (for any choice of length function),
  and we claim that it is thus valid for $H$; this will complete our
  proof.
  
  We may choose a length function
  $\norm\cdot_M$ on $X_*(M)$ prolonging the length function
  $\norm\cdot_H$ on $X_*(H)$ as in Lemma \ref{lem:length-subgroup}.
  Let $X \in \Lie(H)$ be nilpotent; regarding $X$ as an element of
  $\Lie(M)$, we may find $\phi \in X_*(M)$ as in the statement of the
  proposition. According to (3), we have $\phi \in X_*(H)$.  We are going to
  verify that $\phi$ satisfies the conclusion of the proposition for
  $H$ and $X \in \Lie(H)$. Again, property (1) needs no further comment.
  For (2), note first that $M = H\cdot Z$ where $Z$ denotes the center
  of $M$. Then $C_M(X) = C_H(X)\cdot Z$. Setting $C_\phi = C_M(X) \cap
  C_M(\phi(\G_m))$ and $C_\phi' = C_H(X) \cap C_H(\phi(\G_m))$, we
  have $C_\phi = C_\phi' \cdot Z$. Thus the unipotent radical of
  $C_\phi'$ is a normal subgroup of $C_\phi$; since $C_\phi$ is
  reductive, so is $C_\phi'$. This suffices to verify (2).  The
  verification of (3) is similar.  Thus, it only remains to see that
  $\phi$ is optimal for $X$ in $H$. 
  Since $X_*(H) \subset X_*(M)$, optimality of $\phi$ for $H$ follows
  at once.
\end{proof}

\begin{rem}
  Let $\phi$ be as in the proposition.  Since $X \in \glie(2;\phi)$, it
  is clear that either $\phi$ is primitive, or $\frac{1}{2}\phi \in
  X_*(G)$ is primitive (and again optimal for $X$).
\end{rem}

\subsection{Cocharacters associated to nilpotent elements.}

In this subsection, we again suppose that we have fixed a length
function on $X_*(G)$.

Let $X \in \glie$ be nilpotent.  A cocharacter $\phi:\G_m \to G$ is said
to be associated with  $X \in \glie$ if $\Ad(\phi(t))X =
t^2X$ for each $t \in \G_m$, and if $\phi$ takes values in the derived
group of a Levi subgroup $L$ of $G$ for which $X \in \Lie(L)$ is
distinguished.

\begin{prop}
  \label{prop:assoc}
  \begin{enumerate}
  \item There exists a cocharacter which is both optimal for and
    associated with $X$.
  \item If the cocharacter $\phi$ is associated to $X$, then $\Int(g)
    \circ \phi$ is associated to $X$ for each $g \in C_G(X)$.
    Conversely, any two cocharacters associated to $X$ are conjugate
    by $C_G^o(X)$.
  \item If $\phi$ is a cocharacter associated with the nilpotent $X$,
    then the parabolic subgroup $P(\phi)$ coincides with the
    instability parabolic $P(X)$.
  \end{enumerate}
\end{prop}

\begin{proof}
  The optimal cocharacter found by Premet in Proposition
  \ref{prop:premet} is associated with $X$ (by (1) and (3) of that
  proposition). This proves (1).  Assertion (2) follows from
  \cite{jantzen:Nilpotent}*{Lemma 5.3(b)}.  With $\psi$ as in (1) and
  $\phi$ as in (3), (2) implies that $\Int(g) \circ \phi = \psi$ is
  optimal for $X$ for some $g \in C_G^o(X)$.  By
  Proposition\ref{prop:kr}, we have $C_G^o(X) \subset P(X)$. Thus
  $P(X)=P(\psi)=P(\Int(g)\circ\phi) = P(\phi)$, whence (3).
\end{proof}

\begin{rem}
  \label{rem:cochar-existence}
  A proof of the existence of a cocharacter associated with $X$ can be
  extracted from the work by Pommerening (which depends on some
  case-checking for exceptional types); see the overview in
  \cite{jantzen:Nilpotent}*{\S 4}.  The proof given in
  \cite{jantzen:Nilpotent} of part (2) of the proposition is
  elementary: it does not depend on the existence of a cocharacter.
\end{rem}

Write $C=C_G(X)$, let $P=P(X)$ denote the instability parabolic of $X$,
and let $U$ be the unipotent radical of $P$.

\begin{cor}
  \label{cor:Levi-decomposition}
  Let $\phi$ be associated with $X$, and let $R =
  R_u(C)$ be the unipotent radical of $C$.
  \begin{enumerate}
  \item The centralizer $C_\phi$ of $\phi(\G_m)$ in $C$ is reductive,
    and $C = C_\phi \cdot R$.
  \item $R = C \cap U = \{g \in C \mid \lim \phi(t)g\phi(t^{-1}) = 1\}$ and
    $\Lie(R) = \bigoplus_{i \ge 1} \Lie(C)(i;\phi)$.
  \end{enumerate}
\end{cor}

\begin{proof}
  It follows from Premet's result Proposition \ref{prop:premet} that 1
  and 2 are valid for a particular cocharacter $\phi$ associated to
  $X$; the general case results from the conjugacy under $C_G^o(X)$ of
  associated cocharacters.
\end{proof}

\begin{theorem}
  \label{theorem:assoc}
  Let $X\in \glie$ be nilpotent, and let $\phi$ be a cocharacter
  associated to $X$. Then $\phi$ is optimal for $X$.  Conversely,
  suppose that $\psi \in X_*(G)$ is primitive, $\psi$ is optimal for
  $X$, and $X \in \glie(m,\psi)$ for some $m \in \Z_{\ge 1}$. Then
  $m=1$ or $2$.  If $m=2$, then $\psi$ is associated with $X$, if
  $m=1$ then $2\psi$ is associated with $X$.
\end{theorem}

\begin{proof}
  Let $\phi_0$ be a cocharacter which is both optimal for and
  associated with $X$ as in Proposition \ref{prop:assoc}(1).
  
  Suppose first that $\phi$ is associated to $X$. By Proposition
  \ref{prop:assoc}(2), $\phi$ is conjugate under $C_G^o(X)$ to
  $\phi_0$.  Since $C_G^o(X)$ is contained in the instability
  parabolic $P(X)$ by Proposition \ref{prop:kr} (4), optimality of
  $\phi$ follows from Proposition \ref{prop:kr}(3).
  
  Now suppose that $\psi$ is primitive and optimal for $X$, and that
  $X \in \glie(m;\psi)$ as above. Let $P=P(X)$ be the instability
  parabolic, and let $U$ be its unipotent radical. If $\phi_0$ is
  primitive, write $\lambda = \phi_0$. Otherwise we put $\lambda =
  \frac{1}{2}\phi_0$. Thus in each case $\lambda$ is primitive and
  optimal for $X$, and $X \in \glie(n;\lambda)$ with $n=1$ if $\phi_0$
  is not primitive, and $n=2$ if $\phi_0$ is primitive.
  
  By Proposition \ref{prop:kr}(3) $\psi$ and $\lambda$ are conjugate
  via $P$. By \cite{springer-LAG}*{13.4.2}, the centralizer
  $C_G(\psi(\G_m))$ is a Levi subgroup of $P=P(\psi)$. It follows that
  $\psi$ and $\lambda$ are conjugate by an element $u \in U$.  By
  Lemma \ref{lem:parabolic-homog}(2), we see that $m=n$ and $u \in
  C_G(X)$.  Applying Proposition \ref{prop:assoc}(2)
  completes the proof.
\end{proof}

\begin{cor}
  \label{cor:assoc-unicity}
  Let $S$ be a maximal torus of the instability parabolic $P$.  There
  is at most one $\phi \in X_*(S)$ which is associated to $X$.
\end{cor}

\begin{proof}
  Suppose $\phi,\phi' \in X_*(S)$ are associated to $X$.  By the
  previous result, $\phi$ and $\phi'$ are optimal for $X$.  If $\psi$
  denotes the unique primitive optimal cocharacter in $X_*(S)$
  associated with $X$, then $\phi = n\psi$ and $\phi'=n'\psi$ for some
  $n,n' \in \Z_{\ge 1}$. Since $X \in \glie(2;\phi)$ and $X \in
  \glie(2;\phi')$, we see that $n=n'=1$ or 2, and so $\phi=\phi'$.
\end{proof}

Note that while the proof of the preceding corollary  depends on the
choice of the length function on $X_*(G)$, the conclusion
is independent of that choice.

\section{Rationality of associated cocharacters}
\label{sec:cochar-rational}

If $A$ is a linear algebraic group defined over the ground field $F$,
we may always find a maximal torus of $A$ which is defined over $F$;
cf. \cite{springer-LAG}*{13.3.6}. Moreover, any two maximal tori of
$A$ are conjugate by an element of $A^o$ [and even by an element of
$A^o(F_\sep)$]; \cite{springer-LAG}*{Theorem 6.4.1, Prop. 13.3.1}. We
will use these facts without further reference.  In this section,
$G=G_{/F}$ is a reductive group defined over $F$.  We assume
throughout that the characteristic of $k$ is good for $G$.

\subsection{A separability lemma}
\label{sub:separability}

Let $(\rho,V)$ be a linear representation for $G$,
and let $0 \not = v \in V$. Make the following assumptions.
\begin{enumerate}
\item[H1.] Suppose that the
$G$-orbit $\orbit = \rho(G)v$ is separable.
\item[H2.] Suppose that $\orbit$ contains $k^\times w = \{aw \mid a
  \in k^\times\}$ for each $w \in \orbit$.
\end{enumerate}

Observe that H1 and H2 are geometric conditions; they only depend on
$G$ and $V$ over $k$. Recall as well that a $G$-orbit $\orbit$ is
separable just in case some (hence any) orbit map $(g \mapsto
\rho(g)x):G \to \orbit$ for a fixed $x \in \orbit$ has surjective differential
at the identity of $G$.

Denote by $\rho'$ the action
$(g,[w]) \mapsto [\rho(g)w]$ of $G$ on $\Proj(V)$, where $[w]$
denotes the line in $V$ through $w$. Write $\orbit' = \rho'(G)[v]$.

Let $\LL$ be the line bundle over $\Proj(V)$ corresponding to the
invertible coherent sheaf $\mathscr{O}_{\Proj(V)}(1)$.  The
$\G_m$-bundle $\pi:V \setminus \{0\} \to \Proj(V)$  is
the bundle $\LL^\times$ obtained by discarding the zero section from
$\LL$.  In view of H2, $\hat \pi = \pi_{\mid \orbit}:\orbit \to
\orbit'$ is the pull-back of $\LL^\times$ along the inclusion
$i:\orbit' \to \Proj(V)$; in particular it is a (Zariski) locally
trivial principal $\G_m$-bundle.  It follows that $d\hat
\pi_w:T_w \orbit \to T_{\hat \pi(w)} \orbit'$ is surjective for each
$w \in \orbit.$

\begin{lem}
  \label{lem:separability}
  The orbit $\orbit' = \rho'(G)[v]$ is separable. 
\end{lem}

\begin{proof}
  Let $f:G \to \orbit$ be the orbit map $f(g) = \rho(g)v$, and let
  $f':G \to \orbit'$ be $f'(g) = \rho'(g)[v]$.  Since $\orbit$ is
  separable, $df_1:T_1G \to T_v\orbit$ is surjective.  Since $f' =
  \hat \pi \circ f$, and since $d\hat \pi_v$ is surjective, we deduce
  that $df'_1$ is surjective. This proves the lemma.
\end{proof}

\begin{rem}
  The conclusion of the lemma is in general not true when H2 (or
  H1) doesn't hold. Consider the linear representation $(\rho,V)$ of $G
  = \G_m$ where $V=k^2$ and $\rho$ is given by $\rho(t)(a,b) = (t^{-1}a,t^{p-1}b)$.
  Each $G$-orbit on $V$ is separable. However, $\rho'(t)[1:1] =
  [1:t^p]$, hence the (open) orbit $\rho'(G)[1:1]
  \subset \Proj(V)$ is not separable  (the orbit map has 0 differential).
\end{rem}

Consider now the adjoint representation $(\rho,V) = (\Ad,\glie)$ of $G$.
According to \cite{jantzen:Nilpotent}*{\S 2.10, 2.11}, condition H2 of
\S \ref{sub:separability}  is valid for each nilpotent orbit (this holds
even in bad characteristic; the only thing required is the finiteness
of the number of nilpotent orbits over $k$. That finiteness is known
by an uniform argument for good primes, and by case-checking
(Holt--Spaltenstein) for bad primes).

The validity of condition H1 is discussed in Jantzen's notes
\cite{jantzen:Nilpotent}*{\S2.9}. For example, it is valid for the
\emph{standard groups} from \S~\ref{sec:general-reductive}; see
Proposition \ref{prop:standard=>separable}.

\subsection{Associated cocharacters over a ground field.}
\label{sub:rational}

Recall that the characteristic $p$ is assumed to be good for the
reductive $F$-group $G=G_{/F}$.

Fix $X \in \glie(F)$ nilpotent.  We make the following assumption:
\begin{equation}
  \label{eq:F-assumption}
  \text{either $F$ is perfect, or the $G$-orbit of $X$  is separable.}
\end{equation}

Let $N=N(X) = \{g \in G \mid
\Ad(g)X \in kX\}$.  Thus $N$ is the stabilizer of $[X] \in
\Proj(\glie)$ (see \S \ref{sub:separability}).  If $\phi$ is an
cocharacter of $G$ associated with $X$, then $\phi \in X_*(N)$.
Moreover, $\phi(\G_m)$ normalizes $C=C_G(X)$.

\begin{lem}
  \label{lem:N-lemma}
  Let $S$ be any maximal torus of $N$. Then there is a unique
  cocharacter in $X_*(S)$ associated with $X$.
\end{lem}

\begin{proof}
  Fix a cocharacter $\phi$ associated to $X$; then $N = \phi(\G_m)
  \cdot C$ where $C=C_G(X)$; cf. \cite{jantzen:Nilpotent}*{\S5.3}.
  Choose a maximal torus $T$ of $N$ with $\phi(\G_m) \subseteq T$.  If
  $S$ is another maximal torus of $N$, then $S = gTg^{-1}$ with $g \in
  N$. Writing $g^{-1} = \phi(a)h^{-1}$ with $a \in k^\times$ and
  $h \in C_G(X)$, we see that $gTg^{-1} = hTh^{-1}$. It follows that
  $\phi'=\Int(h)\circ \phi$ is a cocharacter of $S$; since $h$
  centralizes $X$, $\phi'$ is associated to $X$ by Proposition
  \ref{prop:assoc}(2).
  
  Since $C_G(X) < P$ by Proposition \ref{prop:kr}, we have $N < P$.
  Thus $S$ is contained in a maximal torus of $P$, and uniqueness of
  $\phi' \in X_*(S)$ then follows from Corollary
  \ref{cor:assoc-unicity}.
\end{proof}

If \eqref{eq:F-assumption} holds, the discussion in
\S\ref{sub:separability} shows that the $G$-orbit of $[X] \in
\Proj(\glie)$ is separable; thus \cite{springer-LAG}*{12.1.2} implies
that the group $N = N(X)$ is defined over $F$.

\begin{theorem}
  \label{theorem:assoc-cochar-over-F}
  Let $X \in \glie(F)$ be nilpotent, and assume that
  \eqref{eq:F-assumption} holds.  Then there is cocharacter $\phi$
  associated to $X$ which is defined over $F$.
\end{theorem}

\begin{proof}
  Since $N$ is defined over $F$, we may choose a maximal torus $S
  \subset N$ defined over $F$.  Let $\phi \in X_*(S)$ be the unique
  cocharacter which is associated to $X$; see Lemma \ref{lem:N-lemma}.
  It follows from \cite{springer-LAG}*{13.1.2} that $\phi$ is defined
  over a separable closure $F_\sep$ of $F$ in $k$. We will show that
  $\phi$ is defined over $F$.
  
  Since $S$ is an $F$-torus, the Galois group $\Gamma$ acts on
  $X_*(S)$: for $\gamma \in \Gamma$, $\psi \in X_*(G)$ and $t \in
  F_\sep$, one has
  \begin{equation}
    \label{eq:gal-on-cochars}
    (\gamma \cdot \psi)(t) = \gamma(\psi(\gamma^{-1}(t))).
  \end{equation}
  We must show that $\psi$ is fixed by each $\gamma \in \Gamma$.  
  To do this, we show that $\gamma \cdot \psi$ is a cocharacter associated to
  $X$.

  First, note that since $X = \gamma(X)$ we have
  \begin{equation*}
    \Ad((\gamma \cdot \psi)(t))X = \gamma(\Ad(\psi(\gamma^{-1}(t)))X)
    = \gamma \left(\gamma^{-1}(t^2) X \right )
    =t^2 X.
  \end{equation*}
  Thu $X \in \glie(2;\gamma \cdot \psi)$, and it just remains to show
  that $\gamma \cdot \psi$ takes values in the derived group of some
  Levi subgroup $M$ of $G$ for which $X \in \Lie(M)$ is distinguished.
  
  Since $\phi$ is itself associated with $X$, there is a Levi subgroup $L$ of
  $G$ such that $X \in \Lie(L)$ is distinguished, and such that
  $\phi(\G_m) < (L,L)$. Let $M = \gamma(L)$. Of course, $M$ is again a
  Levi subgroup. Since $\gamma(X) = X$, we have $X \in \Lie(M)$.  The
  equality $C_{M}(X) = \gamma(C_L(X))$ makes clear that $X$ is
  distinguished in $\Lie(M)$. Moreover, $\gamma(L,L) = (M,M)$, so
  it is clear that $\gamma\cdot\phi(\G_m) < (M,M)$.  
  This completes the proof that $\gamma \cdot \phi$ is associated to $X$.
  
  Since $\gamma \cdot \phi \in X_*(S)$ and since $\phi$ is the unique
  cocharacter in $X_*(S)$ associated with $X$, we deduce $\phi =
  \gamma \cdot \phi$ and the theorem is proved.
\end{proof}

\section{The unipotent radical of a nilpotent centralizer}
\label{sec:split-u-radical}

If $A$ is a linear algebraic $F$-group, recall that the Galois
cohomology set $H^1(F,A)$ is by definition $H^1(\Gamma,A(F_\sep))$
where $F_\sep$ is a separable closure of $F$, and $\Gamma =
\Gal(F_\sep/F)$ is the Galois group.  The basic reference for Galois
cohomology is \cite{SerreGC}; see also \cite{springer-LAG}*{\S12.3}.
The set $H^1(F,A)$ classifies torsors (principal homogeneous spaces)
of $A$ over $F$.  It can be defined as the equivalence classes for a
suitable relation on the set $Z^1(F,A) = Z^1(\Gamma,A(F_\sep))$ of
continuous 1-cocycles with values in $A(F_\sep)$; especially, each
$\alpha \in H^1(F,A)$ may be represented by an $a \in Z^1(F,A)$.  When
$A$ is not Abelian, the set $H^1(F,A)$ is not in general a group, but
it does have a distinguished element -- so it is a ``pointed set'' --
which we sometimes write as $1$.  Thus, the notation $H^1(F,A)=1$
means that this set has one element.

Let $G$ be a reductive $F$-group in good characteristic, and let $X
\in \glie(F)$ be nilpotent.  \emph{Assume throughout this section that
\eqref{eq:F-assumption} holds for $X$.}

We begin by noting:
\begin{prop} 
  \label{prop:nilpotent-rationality}
  The instability parabolic $P(X)$ is defined over $F$.
\end{prop}

\begin{proof}
  By Theorem \ref{theorem:assoc-cochar-over-F}, there is an
  $F$-cocharacter $\phi$ associated with $X$. Since $P(X) = P(\phi)$
  by Proposition \ref{prop:assoc}, $P(X)$ is defined over $F$ by
  Remark \ref{rem:rational-parabolic}.
\end{proof}
This result was previously obtained in
\cite{Mc:sub-principal}*{Theorem 15}.  See
\cite{Ramanan-Ramanathan}*{Theorem 2.3} for a related result.

If $A$ is connected and unipotent (and defined over $F$), one says that $A$ is
$F$-split if there is a sequence of $F$-subgroups $1 = A_n \normal
A_{n-1} \normal \cdots \normal A_2 \normal A_1 = A$ such that each
quotient $A_i/A_{i+1}$ is $F$-isomorphic to the additive group
$\G_{a/F}$.

\begin{theorem}
  \label{theorem:split-u-radical}
  Write $C=C_G(X)$ for the centralizer of $X$, and let $R = R_u(C)$ be
  the unipotent radical.  Then $R$ is defined over $F$ and is an
  $F$-split unipotent group.
\end{theorem}

\begin{proof}
  Let $P=P(X)$ be the instability parabolic of $X$. By Proposition
  \ref{prop:nilpotent-rationality}, $P$ is defined over $F$. Denote by
  $U$ the unipotent radical of $P$; it is defined over $F$ as well
  \cite{springer-LAG}*{13.4.2}. By Corollary
  \ref{cor:Levi-decomposition}(2), the unipotent radical of $C$ is
  $R=C \cap U$, and $\Lie(R) = \Lie(C) \cap \Lie(U)$. Thus, it follows
  from \cite{springer-LAG}*{12.1.5} that $R$ is
  defined over $F$.
  
  By Theorem \ref{theorem:assoc-cochar-over-F}, we may find a
  cocharacter $\phi \in X_*(P)$ associated to $X$ which is defined
  over $F$.  Let $S$ denote the image of $\phi$; it is a 1-dimensional
  split $F$-torus.  It is clear that $S$ acts as a group of
  automorphisms of $R$.  Since $R=U \cap C$, Lemma
  \ref{lem:phi-no-fixed-points} implies that the $F$-torus $S$ has no
  non-trivial fixed points on $R$. It now follows from
  \cite{springer-LAG}*{Corollary 14.4.2} that $R$ is an $F$-split
  unipotent group.
\end{proof}

\begin{cor}
  \label{cor:centralizer-levi-decomp}
  Let $C=C_G(X)$ be the centralizer of $X$. Choose a cocharacter
  $\phi$ associated with $X$ and defined over $F$ (Theorem
  \ref{theorem:assoc-cochar-over-F}).  Then $C = C_\phi \cdot R$ is a
  Levi decomposition defined over $F$, where $C_\phi$ is as in
  Corollary \ref{cor:Levi-decomposition} and $R=R_u(C)$.
\end{cor}

\begin{proof}
  One knows that $C=C_\phi \cdot R$ is a Levi decomposition over $k$;
  the only thing to check is the rationality.  The theorem shows that
  $R$ is defined over $F$. Since $\phi(\G_m)$ is an $F$-torus, its
  centralizer in $C$ is defined over $F$
  (\cite{springer-LAG}*{13.3.1}), whence the corollary.
\end{proof}

\begin{prop}
  \label{prop:split-coh}
  Suppose that $U$ is an $F$-split unipotent group.
  \begin{enumerate}
  \item $H^1(F,U) = 1$, and if $U$ is commutative, $H^i(F,U)=1$ for
    all $i \ge 1$.
  \item If $U$ is a normal subgroup of the $F$-group $A$, and $z \in
    Z^1(F,A)$, then $H^1(F,\twist{z}{U}) =1$, where $\twist{z}{U}$
    denotes the group obtained from $U$ by twisting with $z$.
    If $U$ is commuatative, $H^i(F,\twist{z}{U})=1$ for all $i \ge 1$.
  \end{enumerate}
\end{prop}

\begin{proof}
  Since $U$ has a filtration by normal $F$-subgroups such that each quotient
  is $F$ isomorphic to $\G_{a/F}$, the first assertion follows from the
  additive version of Hilbert 90 \cite{SerreGC}*{II.1.2 Prop. 1} together with a long
  exact sequence argument.
  
  Since $U$ and $\twist{z}{U}$ are isomorphic over $F_\sep$ by construction,
  the second assertion follows from the fact that a unipotent
  $F$-group $V$ is $F$-split if and only if $V_{/F_\sep}$ is
  $F_\sep$-split; see \cite{springer-LAG}*{14.3.8}.
\end{proof}

\begin{rem}
  When $F$ is not perfect, there are $F$-groups $A$ whose unipotent
  radical $R_u(A)$ is not defined over $F$. Take for example the
  $F$-group $A=R_{E/F}\G_m$, where $F \subset E$ is a finite purely
  inseparable extension of degree $p$ and $R_{E/F}$ is Weil's
  restriction of scalars functor. The unipotent radical of $A$ has
  dimension $p-1$, but is not defined over $F$. In fact, $A$ is
  \emph{$F$-reductive}; the maximal closed, connected, normal, unipotent
  $F$-subgroup of $A$ is trivial.
\end{rem}

\begin{rem}
  \label{rem:non-split-example}
  Consider the field $F = \kappa\powfield{t}$ of formal series, where $\kappa$
is any field of characteristic $p > 2$.  Let $U\le \G_a \times \G_a$ be
  the unipotent group $F$-group defined by
  \begin{equation*}
    U = \{(y,z) \in \G_a \times \G_a \mid y^p - y = tz^p\}.
  \end{equation*}
  Then $U$ is defined over $F$, and $U$ is isomorphic over an
  algebraic closure $\overline{F}$ to $\G_{a/\overline{F}}$ (but not
  over $F_\sep$). In fact, $U$ is isomorphic with $\G_a$ over
  $F(t^{1/p})$.  There is an exact sequence
  \begin{equation*}
    F \times F \xrightarrow{(y,z) \mapsto y^p - y - tz^p} 
    F \xrightarrow{\delta} H^1(F,U).
  \end{equation*}
  It is straightforward to verify that the equation $y^p - y = g$ has
  no solution $y \in F$ in case $v(g) < 0$ and $v(g) \not \congruent 0
  \pmod p$, where $v$ denotes the usual $t$-adic valuation on $F$.
  Since $v(tz^p) \congruent 1 \pmod p$ for any $z \in F^\times$, it
  follows that the elements $\{\delta(t^{-np+2}) \mid n \ge 1\}$ of
  $H^1(F,U)$ are all distinct. Thus, $H^1(F,U)$ is infinite; in
  particular, it is non-trivial.  As a consequence, $U$ is not
  $F$-isomorphic to $\G_{a/F}$ and so isn't $F$-split (this is
  \cite{SerreGC}*{II.\S2.1 Exerc.  3}).
\end{rem}

\begin{prop}
  \label{prop:split-unipotent-conj}
  Let the $F$-split unipotent group $U$ act
  on the $F$-variety $X$ (by $F$-morphisms).  Suppose $x,y \in X(F)$
  are conjugate by $U(\overline{F})$.
  Assume:
  \begin{enumerate}
  \item The $U$-orbit of $x$ is separable.
  \item $U_x=\Stab_U(x)$
    is $F$-split.
  \end{enumerate}
  Then $x$ and $y$ are conjugate via $U(F)$.
\end{prop}

\begin{proof}
  Let $\orbit \subset X$ be the orbit $U.x$. Then $\orbit$ is a
  locally closed subvariety of $X$ defined over $F$. Since the orbit
  map $U \to \orbit$ is separable, the group $U_x$ is smooth and there
  is a $U$-equivariant $F$-isomorphism $\orbit \iso U/U_x$.

  Thus there is an exact sequence of pointed sets
  \begin{equation*}
    U_x(F) \to U(F) \to \orbit(F) \to H^1(F,U_x);
  \end{equation*}
  see 
  \cite{springer-LAG}*{12.3.4} (and see the discussion in the beginning
  of \S \ref{sec:finiteness} below).  Since $U_x$ is $F$-split, the
  latter set is trivial. Thus the orbit map $U(F) \to \orbit(F)$ is
  surjective, whence the proposition.
\end{proof}

Recall our assumption \eqref{eq:F-assumption} on the nilpotent element
$X \in \glie(F)$.
\begin{prop}
  \label{prop:arithmetic-U-orbit}
  Let $\phi$ be a cocharacter
  associated with $X$ which is defined over $F$; cf. Theorem
  \ref{theorem:assoc-cochar-over-F}.  Let $\ulie$ be the Lie algebra
  of the unipotent radical $U$ of $P=P(X)$, and let $\vlie =
  \bigoplus_{i \ge 3} \glie(i;\phi) \subset \ulie$.  Then $\Ad(U(F))X
  = X + \vlie(F)$.
\end{prop}

\begin{proof}
  The group $U$ acts on the $F$-variety $X + \vlie$; see the proof of
  Lemma \ref{lem:parabolic-homog}. Moreover, the stabilizer $U_X =
  C_U(X)$ is precisely the unipotent radical of $C_G(X)$; see
  Corollary \ref{cor:Levi-decomposition}.  Especially, the $U$-orbit
  of $X$ is separable, and $U_X$ is $F$-split. Thus the result will
  follow from the previous proposition provided that $\Ad(U)X = X + \vlie;$
  i.e. that the proposition holds over $k$.
  
  Well, by \cite{jantzen:Nilpotent}*{Prop. 5.9(c)}, we have
  $\overline{\Ad(P)X} = \bigoplus_{i \ge 2} \glie(2;\phi)$.  Since the
  orbit of $X$ is separable, the differential of the orbit map is
  surjective.  Thus
  \begin{equation}
    \label{eq:ad-onto}
    \ad(X):\glie(0) \to \glie(2) \quad \text{and} \quad 
    \ad(X):\ulie \to \vlie \quad \text{are surjective}.
  \end{equation}
  From \eqref{eq:ad-onto}, we see that the orbit map $\rho:U \to X
  + \vlie$ given by $u \mapsto \Ad(u)X$ is dominant. By a result of
  Rosenlicht, one knows that each $U$-orbit on
  the affine variety $X + \vlie$ is closed (see
  \cite{MR50:4766}*{Prop. 2.5}). Since $\Ad(U)X$ is dense
  in the irreducible variety $X + \vlie$, equality follows.
\end{proof}

\section{Galois cohomology and finiteness}
\label{sec:finiteness}

Let $X \in \glie(F)$ and suppose that \eqref{eq:F-assumption} holds.
Since the centralizer $C=C_G(X)$ is smooth, there is a $G$-equivariant
$F$-isomorphism $\orbit = \Ad(G)X \iso G/C$.  We have thus an exact
sequence of pointed sets
\begin{equation}
  \label{eq:orbit-seq}
  C(F) \to G(F) \to \orbit(F) \xrightarrow{\delta} H^1(F,C) 
  \xrightarrow{\alpha} H^1(F,G);
\end{equation}
see \cite{springer-LAG}*{Prop. 12.3.4}.
One should be cautious concerning such exact sequences.  A sequence
of pointed sets $A \xrightarrow{f} B \to 1$ is exact if and only if
$f$ is surjective.  On the other hand, the sequence of pointed sets $1
\to A \xrightarrow{f} B$ can be exact even when $f$ is not injective.
Using techniques of ``twisting'' one shows that the $G(F)$-orbits in
$\orbit(F)$ are in bijection with the kernel of $\alpha:H^1(F,C) \to
H^1(F,G)$; see \cite{b-invol}*{28.2} or \cite{SerreGC}*{\S5.4 Cor. 2}.

\begin{lem}
  \label{lem:split-radical-finiteness}
  Let $A$ be a linear algebraic $F$-group, and suppose that $R$ is a
  normal connected unipotent $F$-subgroup which is $F$-split.
  \begin{enumerate}
  \item If $H^1(F,A/R)$ is finite, then $H^1(F,A)$ is finite.
  \item Suppose that an $F$-split torus $S$ acts on $A$ as a group of
    automorphisms, and assume that 1 is the only fixed point of $S$ on
    $R$. Then the natural map $H^1(F,A) \to H^1(F,A/R)$ is bijective.
  \end{enumerate}
\end{lem}

\begin{proof}
  Since $R$ is $F$-split, $H^1(F,\twist{z}{R}) = 1$ for each $z \in
  Z^1(F,A)$ by Proposition \ref{prop:split-coh}. We have an exact
  sequence $1 \to H^1(F,A) \xrightarrow{f} H^1(F,A/R)$. Fix $\alpha
  \in H^1(F,A)$. By \cite{SerreGC}*{I.5 Cor. 2} the elements $\beta \in
  H^1(F,A)$ with $f(\alpha)=f(\beta)$ are in bijection with a certain
  quotient of $H^1(F,\twist{a}{R})$, where $a \in Z^1(F,A)$ represents
  $\alpha$.  Thus $f$ is injective and (1) follows.
  
  Supposing now that a split torus $S$ acts as in (2), we show that
  $f$ is surjective. Let $R_1 = (R,R)$ be the derived group, and for
  $i>1$, let $R_i = (R,R_{i-1})$.  Then $R_n =\{1\}$ for some $n \ge
  1$, and each $R_i$ is normal in $A$. In particular, $S$ acts without
  non-trivial fixed points on each $R_i$, so that each $R_i$ is an
  $F$-split unipotent group by \cite{springer-LAG}*{14.4.2}. Moreover,
  each $R_i/R_{i+1}$ is an $F$-split commutative unipotent group, by
  \cite{springer-LAG}*{14.3.12 exercise 2}.
  
  If we show that the natural map $H^1(F,A) \xrightarrow{f}
  H^1(F,A/R_{n-1})$ is a bijection, the result for $R$ will 
  follow by induction. Thus, we may suppose  $R$ to be $F$-split and
  \emph{commutative}, and we must show that $f$ is surjective.  In
  this case, we have $H^2(F,\twist{z}{R})=1$ for all $z \in
  Z^1(F,A/R)$ by Proposition \ref{prop:split-coh}.  Thus, we may apply
  \cite{SerreGC}*{I.5 Cor. to Prop 41} to see that $f$ is surjective.
\end{proof}

\begin{lem}
  \label{lem:reductive-cohomology}
  Suppose that $F$ has cohomological dimension $\le 1$. 
  Let $A$ be a linear algebraic $F$-group, and suppose that the
  $F$-group $A^o$ is reductive. Then the natural map
  \begin{equation*}
    H^1(F,A) \xrightarrow{f} H^1(F,A/A^o)
  \end{equation*}
  is injective. The map $f$ is bijective if moreover $F$ is perfect.
\end{lem}

\begin{proof}
  Since $F$ has cohomological dimension $\le 1$, a result of
  Borel and Springer \cite{BS68}*{8.6} implies tht $H^1(F,A^o) = 0$. Thus
  the exact sequence in Galois cohomology arising from the sequence
  \begin{equation*}
    1 \to A^o \to A \to A/A^o \to 1
  \end{equation*}
  shows that $1 \to H^1(F,A) \xrightarrow{f} H^1(F,A/A^o)$ is exact.
  The proof that $f$ is injective may then be found in the proof of
  \cite{SerreGC}*{III.2.4 Corollary 3}; note that $F$ is assumed
  perfect in \emph{loc. cit.} but this is not essential for the proof
  of injectivity [one just needs to use: if $b \in Z^1(F,A)$, the
  $F$-group $\twist{b}{A^o}$ is again connected and reductive and
  hence has trivial $H^1$ by the result of Borel--Springer]. This same
  result shows that $f$ is bijective in case $F$ is also perfect.
%  Fix $\beta \in H^1(F,A)$ and let $b \in Z^1(F,A)$ represent $\beta$.
%  According to \cite{SerreGC}*{I.5.6 cor. 2}, the elements $\alpha \in
%  H^1(F,A)$ with $f(\alpha) = f(\beta)$ correspond bijectively with
%  elements in a certain quotient of $H=H^1(F,\twist{b}{(A^o)})$, where
%  $\twist{b}{(A^o)}$ denotes the group obtained from $A^o$ by twisting
%  with the cocycle $b$.  Since $\twist{b}{(A^o)}$ is again a connected
%  and reductive $F$-group, $H=1$ by the result of Borel--Springer
%  cited before, and so $f(\alpha) = f(\beta)$ implies $\alpha =
%  \beta$. Thus $f$ is indeed injective.
\end{proof}

In the previous proof, the surjectivity of $f$ when $F$ is perfect
depends on a result of Springer \cite{SerreGC}*{III.2.4 Theorem 3}
concerning principal homogeneous spaces.

Recall that we suppose \eqref{eq:F-assumption} to hold for the
nilpotent $X \in \glie(F)$.  The centralizer $C$ is then defined over
$F$, and hence the connected component $C^o$ of $C$ is defined over
$F$ as well. We write $A_X$ for the component group $C/C^o$; it is a
finite linear $F$-group. 

\begin{prop}
  Suppose that $F$ has cohomological dimension $\le 1$.  Then:
  \begin{enumerate}
  \item Each element of $A_X(F)$ can be represented by a coset $gC^o$
    with $g \in C(F)$.
  \item The set of $G(F)$-orbits in $\orbit(F)$
    identifies with a subset of $H^1(F,A_X)$.
  \end{enumerate}
\end{prop}

\begin{proof}
  Let $\phi$ be a cocharacter associated to $X$ which is defined over
  $F$. Then the $F$-split torus $S=\phi(\G_m)$ acts as a group of
  automorphisms of $C$ and the only fixed point on the unipotent
  radical $R$ of $C$ is the identity. So Lemma
  \ref{lem:split-radical-finiteness} shows that the natural maps
  $H^1(F,C) \xrightarrow{g} H^1(F,C/R)$ and $H^1(F,C^o) \to
  H^1(F,C^o/R)$ are bijective.
  
  Now, $C^o/R$ is a connected, reductive $F$-group, so $H^1(F,C^o/R) =
  1$ by the result of Borel--Springer \cite{BS68}*{8.6} cited in the
  previous Lemma.  It follows that $H^1(F,C^o) = 1$.  There is thus an
  exact sequence
  \begin{equation*}
    1 \to C^o(F) \to C(F) \to A_X(F) \to 1 
  \end{equation*}
  which proves (1).  
  
  It follows from \cite{springer-steinberg}*{III.3.15} that each class
  in $A_X$ can be represented by a semisimple element; thus $A_X \iso
  (C/R)/(C^o/R)$.  
  It follows from Lemma \ref{lem:reductive-cohomology} that the map
  \begin{equation*}
    H^1(F,C/R) \xrightarrow{f} H^1(F,A_X)
  \end{equation*}
  is injective.
  
  If $f':H^1(F,C) \to H^1(F,A_X)$ is induced by the quotient map $C
  \to A_X$, then $f'=f \circ g$. Since $g$ is bijective, $f'$ is
  injective.  After the result of Borel--Springer already cited, we
  have $H^1(F,G)=1$; thus \cite{b-invol}*{Cor 28.2} implies that the
  map $\delta:\orbit(F) \to H^1(F,C)$ from \eqref{eq:orbit-seq} induces a
  bijection between $H^1(F,C)$ and the set of $G(F)$-orbits in $\orbit(F)$.
  Assertion (2) now follows.
\end{proof}

\begin{rem}
  With assumptions and notation as in the preceding proposition, if
  $F$ is perfect one knows by Lemma \ref{lem:reductive-cohomology}
  that the set of $G(F)$-orbits in $\orbit(F)$ identifies with
  $H^1(F,A_X)$. One might well wonder if this remains so when $F$ is
  not assumed perfect (assuming \eqref{eq:F-assumption} to hold, of course).
\end{rem}

Let $F$ be a field complete with respect to a non-trivial discrete valuation
$v$, with $p =\char\ F$.  By the residue field $\kappa$ we mean the
quotient of the ring of integers of $F$ by its unique maximal ideal.

\begin{prop}
  \label{prop:finiteness}
  Suppose $\kappa$ is finite or algebraically closed, and let $A$ be a
  linear algebraic group over $F$. Suppose further that
  \begin{enumerate}
  \item The unipotent radical of $A$ is defined
  over $F$ and is $F$-split.
  \item $|A/A^o|$ is invertible in $F$.
  \end{enumerate}
  Then $H^1(F,A)$ is finite.
\end{prop}

\begin{proof}
  Let $n=|A/A^o|$.  Since $n$ is invertible in $F$, the field $F$ has
  only finitely many extensions of degree $n$ in a fixed separable
  closure $F_\sep$. Indeed, there is $\le 1$ unramified extension $F
  \subset F_m$ of degree $m$ for each $m|n$. So we just need to show
  that the number of totally ramified extensions of $F_m \subset F'$
  of degree $n/m$ is finite.  When the residue field is finite, this
  follows from Krasner's Lemma \cite{serre79:_local_field}*{II.2 exer.
    1,2} (since $n/m$ is prime to $p$, the space of separable
  Eisenstein polynomials of degree $n/m$ is compact).  When the
  residue field is algebraically closed, apply
  \cite{serre79:_local_field}*{IV, \S2, Cor. 1, Cor. 3}.
  
  It follows that $H^1(F,A)$ is finite whenever $A$ is a finite
  $\Gamma=\Gal(F_\sep/F)$-group whose order is prime to $\char F$.
  This is a variant of \cite{SerreGC}*{II.4.1 Prop.  8}, and a proof is
  given in \cite{Morris-rational-conjugacy}*{Lemma 3.11}. For
  completeness, we outline the argument here. We may find an open
  normal subgroup $\Gamma_o \normal \Gamma$ which acts trivially on
  $A$ (we suppose the profinite group $\Gamma$ to act continuously on
  $A$).  The open subgroups of $\Gamma_o$ with index dividing $n$ are
  finite in number (apply the conclusion of the previous paragraph to
  the fixed field $F'=(F_\sep)^{\Gamma_o}$), and their intersection is
  an open normal subgroup $\Gamma_1$ of $\Gamma$. Every continuous
  homomorphism $\Gamma_o \to A$ vanishes on $\Gamma_1$, so the
  restriction map $H^1(\Gamma,A) \to H^1(\Gamma_1,A)$ is trivial.  It
  then follows from \cite{SerreGC}*{\S5.8 a)} that $H^1(F,A)$
  identifies with $H^1(\Gamma/\Gamma_1,A)$, which is clearly finite.
  
  Now suppose the proposition is proved in case $A$ is connected.
  Since $\twist{z}{A^o}$ is connected for all $z \in Z^1(F,A)$, and
  since $H^1(F,A/A^o)$ is finite by the previous paragraph, we may
  apply \cite{SerreGC}*{I.5 Cor. 3} to the exact sequence $H^1(F,A^o)
  \to H^1(F,A) \to H^1(F,A/A^o)$ and deduce the proposition for
  general $A$.  Since the unipotent radical of $A$ is split, Lemma
  \ref{lem:split-radical-finiteness} shows moreover that we may
  suppose $A$ to be connected and reductive.
  
  If $\kappa$ is algebraically closed, then by a result of Lang, $F$
  is a $(C_1)$ field; see \cite{SerreGC}*{II.3.3(c)}. In particular,
  $F$ has cohomological dimension $\le 1$; c.f. II.3.2 of \emph{loc.
    cit.}  So when $A$ is connected and reductive, we have
  $H^1(F,A)=1$ by the result of Borel--Springer cited in the proof of
  the previous proposition.
  
  When $\kappa$ is finite, the finiteness of $H^1(F,A)$ for $A$
  connected and reductive is a consequence of Bruhat-Tits theory;
  cf. \cite{SerreGC}*{III.4.3 Remark(2)}.
\end{proof}

\begin{theorem}
  \label{theorem:finiteness}
  Suppose that $F$ is complete for a non-trivial discrete valuation,
  and that the residue field $\kappa$ of $F$ is finite or
  algebraically closed.  If \eqref{eq:F-assumption} holds for the
  nilpotent element $X \in \glie(F)$, then $G(F)$ has finitely many
  orbits on $\orbit(F)$.  In particular, if \eqref{eq:F-assumption}
  holds for each nilpotent $X \in \glie$, the nilpotent $G(F)$ orbits
  on $\glie(F)$ are finite in number.
\end{theorem}

\begin{proof}
  The Bala-Carter-Pommerening theorem implies that there are finitely
  many geometric nilpotent orbits; see \cite{jantzen:Nilpotent}*{\S4}.
  So the final assertion follows from the first.  Now,
  \eqref{eq:orbit-seq} shows that the first assertion follows once we
  know that $H^1(F,C)$ is finite, where $C=C_G(X)$.  The order of the
  component group $A_X = C/C^o$ is invertible in $F$
  \cite{springer-steinberg}*{3.19} (this could also be deduced from the
  explicit results in \cite{sommers-mcninch}).  According to Theorem
  \ref{theorem:split-u-radical}, the unipotent radical of $C$ is
  defined over $F$ and is $F$-split. Thus the theorem follows from the
  previous proposition.
\end{proof}

\begin{rem}
  \begin{enumerate}
  \item Theorem \ref{theorem:finiteness} was obtained by Morris
    \cite{Morris-rational-conjugacy}, under the assumption $p > 4h-4$
    where $h$ denotes the Coxeter number of $G$. The main new
    contribution of the present work is application of Theorem
    \ref{theorem:split-u-radical}.
  \item Recall that \eqref{eq:F-assumption} holds for each nilpotent
    $X \in \glie$ in case $G$ is a \emph{standard} reductive $F$-group;
    cf. Proposition  \ref{prop:standard=>separable}.
  \end{enumerate}
\end{rem}

\section{An example: a non-quasisplit group of type $C_2$}

In this section, we use Proposition \ref{prop:nilpotent-rationality}
to study the arithmetic nilpotent orbits of a group of type $C_2$
which is not quasisplit over the ground field $F$ (i.e. has no
Borel subgroup defined over $F$). In case $F$ is a local field of odd
characteristic, we use some local class field theory to classify these
orbits; we see especially that they are finite in number, as promised
by Theorem \ref{theorem:finiteness}.

Let $Q$ be a division algebra with center $F$ and $\dim_F Q=4$ (one
says that $Q$ is a quaternion division algebra over $F$), and suppose
that $\char F \not = 2$.  There is a uniquely determined symplectic
involution $\iota$ on $Q$; see for example \cite{b-invol}*{\S I.2.C}.

Denote by $A=\Mat_2(Q)$, and let $\sigma$ be the involution
of $A$ given by 
\begin{equation*}
\sigma
\begin{pmatrix}
  \alpha & \beta \\
  \gamma & \delta
\end{pmatrix} = 
\begin{pmatrix}
  \iota(\delta) & \iota(\beta) \\
  \iota(\gamma) & \iota(\alpha)
\end{pmatrix}. 
\end{equation*}
Then $\sigma$ is the adjoint involution determined by an isotropic
hermitian form on a 2 dimensional $Q$-vector space; cf.
\cite{b-invol}*{I.4.A}. 

The algebra $A$ together with the symplectic involution $\sigma$
determine an $F$-form $G_{/F} = \Iso(A,\sigma)$ of $\SP_4$;
we have 
\begin{equation*}
  G(\Lambda) =
  \{g \in A \tensor_F \Lambda \mid g\cdot \sigma(g) =1 \}
\end{equation*}for each commutative $F$-algebra $\Lambda$.
The group $G$ has no Borel subgroup over $F$ (see
\cite{springer-LAG}*{17.2.10}).  There is a
cocharacter $\phi = (t \mapsto
\begin{pmatrix}
  t & 0 \\
  0 & t^{-1}
\end{pmatrix})$
defined over $F$, and $P = P(\phi)$ is a minimal
$F$-parabolic subgroup. By \cite{springer-LAG}*{Theorem 15.4.6}, and a
little thought, all proper $F$-parabolic subgroups of $G$ are
conjugate by $G(F)$.

There are four geometric nilpotent orbits in $\lie{sp}_4(F_\sep)$; the
corresponding conjugacy classes of instability parabolics are all
distinct. So applying Proposition \ref{prop:nilpotent-rationality}, we
see that there is a unique non-0 geometric nilpotent orbit with an
$F$-rational point.

For $0 \not = a \in \Skew(Q,\iota) = \Skew(Q) = \{x \in Q \mid x +
\iota(x)=0\}$, the element
\begin{equation*} X_a =
  \begin{pmatrix}
    0 & a \\
    0 & 0
  \end{pmatrix} \in \glie(F) =\Skew(A,\sigma)
\end{equation*}
is nilpotent.  If the field $L$ splits $Q$, $X_a$ has rank 2 in
$\Mat_4(L) = A \tensor_F L$. It follows from the description of
nilpotent orbits in $\lie{sp}_4$ by partition that $X_a$ lies in the
\emph{subregular orbit} $\orbit_{\text{sr}} = \orbit$ (i.e. $X_a$ acts
with partition $(2,2)$ on the natural symplectic module).

The preceding discussion shows that $\orbit$ is defined over
$F$ and has an $F$-rational point. Moreover, $\orbit(F)$ is precisely
the set of nilpotent elements in $\glie(F)$.

Denote by $M$ the subgroup $\left\{\begin{pmatrix}
    x & 0 \\
    0 & \iota(x^{-1})
\end{pmatrix}\mid x \in \GL_{1,Q}\right\} < P$. Thus
$M \iso \GL_{1,Q}$, and since $M$ is the centralizer of the image of
the cocharacter $\phi$, it is a Levi factor in $P$.  Since a
subregular nilpotent element lies in the Richardson orbit of its
instability parabolic, it follows that \emph{the arithmetic nilpotent
  orbits of $G(F)$ are in bijection with the $M(F)$ orbits on the
  nilradical of $\Lie(P)$}; by the nilradical we mean the Lie algebra
of the unipotent radical of $P$.  Moreover, $P$ is the instability
parabolic for each of the nilpotent elements $X_a$ with $0 \not = a
\in \Skew(Q,\iota)$, and one can even see that $\phi$ is a cocharacter
associated with $X_a$.

The action of $M(F)$ on the nilradical of $\Lie(P)(F)$ identifies with the
representation $(\rho,\Skew(Q))$ of $Q^\times = \GL_{1,Q}(F)$ given by
\begin{equation*}
  \rho(x)y = xy\iota(x) = \dfrac{1}{\Nrd(x)}xyx^{-1},
\end{equation*}
where $\Nrd:Q^\times \to F^\times$ is the reduced norm.
So we seek a description of the $Q^\times$-orbits on
$\Skew(Q)$. One easily sees that the function
\begin{equation*}
  \eta=(y \mapsto \Nrd(y)F^{\times 2}):\Skew(Q)^\times \to F^\times/F^{\times 2}
\end{equation*}
is constant on $Q^\times$-orbits, so the essential problem is to find
the $Q^\times$-orbits on the fibers of $\eta$. If $0 \not = y \in \Skew(Q)$,
$F[y]$ is a maximal subfield of $Q$. It follows that $\eta(y)$ is not
the trivial square class.

Let now $F$ be the local field $\F_q\powfield{t}$ where $\F_q$ is the
finite field having $q$ elements, where $q$ is odd. Then there is a
unique division quaternion algebra $Q$ over $F$
\cite{serre79:_local_field}*{ch.  XIII}.  The group of square classes
$F^\times/F^{\times 2}$ is $\F_q^\times/\F_q^{\times 2} \times \Z/2
\iso \Z/2 \times \Z/2$; cf.  \emph{loc. cit.}  ch. V, Lemma 2.

We claim that the image of $\eta$ consists in the non-trivial elements of
$F^\times/F^{\times 2}$, and that each fiber of $\eta$ is a single
$Q^\times$-orbit. It will follow that there are 3 orbits of $Q^\times$
on $\Skew(Q)^\times$, and thus 3 non-zero arithmetic nilpotent $G(F)$-orbits on
$\glie(F)$.

Let $x \in F^\times$ represent a non-trivial square class, and let
$\Theta < F^\times$ be the subgroup generated by $F^{\times 2}$ and
$x$. Then $\Theta$ is a closed subgroup of index 2 in $F^{\times 2}$
and so $\Theta = N_{L/F}(L^\times)$ for some quadratic extension $L$
of $F$, by local class field theory \cite{serre79:_local_field}*{Ch.
  XIV \S6.  Theorem 1}.  We may write $L=F[\sqrt{a}]$ for some $a \in
F^\times$ with $\Nrd(a) \in xF^{\times2}$. By
\cite{serre79:_local_field}*{Ch.  XIII \S4.  Cor.  3}, $L$ embeds in
$Q$ as a maximal subfield.  Under any such embedding, $\sqrt{a}$
corresponds to an element $y \in \Skew(Q)$ with $\eta(y) =
xF^{\times2}$.  This proves our first claim.

If $y_1,y_2 \in \Skew(Q)^\times$ and $\eta(y_1)=\eta(y_2)$, we show
that $y_1$ and $y_2$ are conjugate under $Q^\times$.  One knows
$F[y_1] \iso F[y_2]$ (as $F$-algebras), so we may suppose, by the
Skolem-Noether Theorem, that $y_1 \in F[y_2]=L$.  Since $y_1,y_2$ each
have trace 0, $y_1\cdot y_2^{-1} \in F^\times$. If $y_1 \cdot y_2^{-1}
= N_{L/F}(\beta)$ for some $\beta \in L^\times$, then $y_1 =
\rho(\beta)y_2$ and our claim holds. If $y_1 \cdot y_2^{-1} \not \in
N_{L/F}(L^\times)$, apply the Skolem-Noether theorem to find $\gamma
\in Q^\times$ such that $y \mapsto \gamma y \gamma^{-1}$ is the
non-trivial element of $\Gal(L/F)$. Then $\{1,\gamma\}$ is an
$L$-basis of $Q$, and moreover, $\gamma^2 \in Z(Q)=F$, so $Q$
identifies with the ``cyclic $F$-algebra'' $(L,\gamma^2)$
\cite{b-invol}*{\S 30.A}. Since $Q$ is a division algebra, $\gamma^2
\not \in N_{L/F}(L^\times)$; see \cite{b-invol}*{Prop.  30.6}. Since
$\Trd(\gamma)=0$, we find
\begin{equation*}
  \rho(\gamma) y_2 = \dfrac{-1}{\Nrd(\gamma)} y_2 = \dfrac{1}{\gamma^2} y_2.
\end{equation*}
Since $[F^\times:N_{L/F}(F^\times)] = 2$, this implies $y_1 \cdot
\rho(\gamma)y_2^{-1} = \dfrac{1}{\gamma^2}y_1y_2^{-1}\in
N_{L/F}(L^\times)$ and the claim follows.

\begin{rem}
  If $0 \not = a \in \Skew(Q)$, the connected component of 1 in the
  centralizer $C=C_G(X_a)$ has dimension 1 and is isomorphic to the
  norm torus $\G^1_{m/L} = \ker(N_{L/F}:R_{L/F}\G_m \to \G_m)$, where
  $L=F[a]$. Moreover, $[C:C^o] = 2$, $C$ is non-abelian, and the
  non-trivial coset of $C^o$ in $C$ has no $F$-rational point.  The
  above calculation shows that $|H^1(F,C)|=3$ when
  $F=\F_q\powfield{t}$.
\end{rem}

\section{Orbital integrals}
\label{sec:orbital-int}

We now suppose that our field $F$ is complete with respect to a
non-trivial discrete valuation $v$, and that the residue field $\f$ is
finite. We suppose that the valuation satifies $v(t)=1$ for a prime
element $t \in F$; the normalized absolute value of $0 \not = a \in F$
is then the rational number $|a| = |\f|^{-v(a)}$.  If the
characteristic of $F$ is 0, we will have nothing new to say in this
section.  When $F$ has characteristic $p>0$, it is isomorphic to the
field of formal power series $\f\powfield{t}$.

If $X$ is a smooth quasi-projective variety over $F$, then $X(F)$ is
an analytic $F$-manifold.
If $\omega$ is a non-vanishing regular differential form
on $X$ of top degree defined over $F$, it defines a measure $|\omega|$
on the locally compact topological space $X(F)$ in a well-known
manner; see e.g. \cite{Platonov}*{\S 3.5}.

Throughout this section, let $G$ be a reductive group defined over
$F$, and suppose that $G$ is $F$-standard.  Recall that
all adjoint orbits and all conjugacy classes are thus known
to be separable; cf. Proposition \ref{prop:standard=>separable}.

Since $G$ is reductive,
the representation of $G$ on $\bwedge^{\dim G} \glie$ is trivial (the
restriction of this representation to a maximal torus of $G$ is
evidently trivial). Thus a left $G$-invariant differential form
$\omega_G$ on $G$ of top degree is also right invariant, so it defines
a left- and right- Haar measure $\abs{\omega_G}$ on the locally
compact group $G(F)$.

Let $X \in \glie(F)$ or $x \in G(F)$, and let $\orbit$ be the
geometric orbit of this element (thus $\orbit \subset \glie$ or
$\orbit \subset G$), and let $C$ be its centralizer.  Since $G$ is
$F$-standard, Proposition \ref{prop:standard=>separable} shows that
$C$ is defined over $F$.

\begin{lem}
  \label{lem:unimodular}
  There is a non-vanishing differential form $\tau$ of top degree on
  $\orbit \iso G/C$ which is defined over $F$. Thus, $C(F)$ is
  unimodular.
\end{lem}

\begin{proof}
  Since there is a $G$-invariant bilinear form on $\glie$ defined over
  $F$, the lemma follows from \cite{springer-steinberg}*{3.24, 3.27}.
\end{proof}

Write $\WW = \Ad(G(F))X$ when $X \in \glie(F)$, and write 
$\WW = \Int(G(F))x$ when $x \in G(F)$.
\begin{lem}
  \label{lem:G-orbit}
  $\WW$ is an open submanifold of $\orbit(F)$, and is a locally closed
  subspace of $\glie(F)$ or of $G(F)$.
\end{lem}

\begin{proof}
  $\orbit$ is a smooth variety defined over $F$, so $\orbit(F)$ is an
  analytic $F$-manifold. We have supposed that the orbit map $G \to \orbit$ is
  separable; in other words, this map has surjective differential
  at each $g \in G$.  The inverse function theorem
  \cite{Serre-LieGroups}*{LG3.9} implies that $\WW$ is open in $\orbit(F)$,
  whence the first assertion.
  
  Now, $\orbit$ is Zariski-open in $\overline{\orbit}$, so $\orbit(F)$
  is open in $\overline{\orbit}(F)$ in the $F$-topology.  Thus $\WW$
  is open in $\overline{\orbit}(F)$, which shows that $\WW$ is open in
  its closure $\overline{\WW} \subset \overline{\orbit}(F)$.
\end{proof}

For a topological space $\X$ we will write $C(\X)$ for the algebra
of $\C$-valued continuous functions on $\X$, and $\cc(\X)$ for the
sub-algebra of compactly supported continuous functions.  

With $\tau$ as in lemma \ref{lem:unimodular}, we obtain a
$G(F)$-invariant measure $dg^*$ on $G(F)/C(F)$.  For $f \in
\cc(\glie(F))$, respectively $f \in \cc(G(F))$, define the
\emph{orbital integral} of $f$ over $\WW$ to be
\begin{equation*}
  I_X(f) = \int_{G(F)/C(F)} f(\Ad(g)Y)dg^* \text{for $X \in \glie(F)$, and}
\end{equation*}
\begin{equation*}
  I_x(f) = \int_{G(F)/C(F)} f(gxg^{-1})dg^* \text{for $x \in G(F)$.}
\end{equation*}
By construction, of course, we have $I_X(f), I_x(f) < \infty$ if $f_{\mid \WW}
\in \cc(\WW)$; this is so e.g. if $\WW$ is \emph{closed}.  One is interested
in the convergence of the integrals $I_X(f)$ in general; we will now investigate
these integrals.

\subsection{Nilpotent case}

We first consider the integral $I_X(f)$ in the case where $X \in
\glie(F)$ is nilpotent.
\begin{theorem}
  \label{theorem:nilpotent-convergence}
  Let $X \in \glie(F)$ be nilpotent. Then 
  $I_X(f) < \infty$ for each $f \in \cc(\glie(F))$.
\end{theorem}

The theorem was proved by Deligne and by Ranga Rao \cite{rao-orbital},
in the case that $F$ has characteristic 0.  We show here how to adapt
the original proof to the positive characteristic setting.

Let $P$ be the instability $F$-parabolic subgroup determined by $X$.
Fix a co-character $\phi$ associated to $X$ and defined over $F$; cf.
Theorem \ref{theorem:assoc-cochar-over-F}.  We abbreviate
$\glie(i;\phi)$ as $\glie(i)$ for $i \in \Z$, and we write $\wlie_i =
\glie(i)(F)$.  Recall that $\phi$ determines a Levi factor $M =
C_G(\phi(\G_m))$ of $P$ which is defined over $F$.

Inspecting the argument given in \cite{rao-orbital}, one sees that
$I_X(f) < \infty$ for $f \in \cc(\glie(F))$ if we establish the
following:
\begin{enumerate}
\item[R1.] The $M(F)$-orbit of $X$ is an open submanifold $\VV \subset \wlie_2$.
\item[R2.] The $P(F)$-orbit of $X$ is $\VV + \sum_{i \ge 3} \wlie_i$.
\item[R3.] There is a non-negative function $\phi \in C(\wlie_2)$
  with $\phi(X) \not = 0$ and $\phi(\Ad(m)Y) =
  \abs{\det(\Ad(m)_{\wlie_1})}\phi(Y)$.
\end{enumerate}

More precisely, suppose that R1--3 hold, let $K$ be an open compact
subgroup of $G(F)$ with the property $G(F) = K \cdot P(F)$ (that there
should be such a $K$ is a result of Bruhat--Tits; see e.g.
\cite{tits:reductive/local}), let $dY$ and $dZ$ be additive Haar
measure respectively on $\wlie_2$ and $\wlie_{\ge 3} = \sum_{i \ge
  3} \wlie_i$, and put
\begin{equation*}
  \Lambda(f) = \int_{\wlie_2 \oplus \wlie_{\ge 3}} 
  \phi(Y)f(Y+Z) dY dZ,
\end{equation*}
and
\begin{equation*}
  \overline{f}(Y) = \int_K f(\Ad(x)Y)dx, \quad Y \in \glie(F),
\end{equation*}
for $f \in \cc(\glie(F))$,
$dx$ denoting a Haar measure on $K$.
Under our assumptions, it is proved in \emph{loc. cit.} that
\begin{equation*}
  I_X(f) = c \cdot \Lambda(\overline{f}) \quad \text{for } 
  f \in \cc(\glie(F)),
\end{equation*}
where $0 \not = c$ is a suitable constant; in particular, $I_X(f) < \infty$.

We first verify that conditions R1, R2 hold.
\begin{prop}
  \label{prop:P-orbit} Let $\vlie_+ = \sum_{i \ge 3} \wlie_i.$
  \begin{enumerate}
  \item The $M(F)$-orbit of $X$ is an open submanifold
    $\VV \subset \wlie_2$. 
  \item The $P(F)$-orbit of $X$ is $\VV + \vlie_+$, an open
    submanifold of $\vlie$.
  \end{enumerate}
\end{prop}

\begin{proof}  
  The orbit map $(m \mapsto \Ad(m)X):M \to \glie(2)$ has surjective
  differential at 1 by \eqref{eq:ad-onto} in the proof of Proposition
  \ref{prop:arithmetic-U-orbit}, and hence has surjective differential
  at each $m \in M$ (by transport of structure).  It follows as in the
  proof of Lemma \ref{lem:G-orbit} that $\VV$ is open.
  
  Now $P(F) = M(F)\cdot U(F)$; moreover, $\Ad(U(F))X = X + \vlie_+$ by
  Proposition \ref{prop:arithmetic-U-orbit}. Thus $\Ad(P(F))X =
  \Ad(M(F))(X + \vlie_+) = \VV + \vlie_+$ as claimed.
\end{proof}
The proposition implies R1 and R2. Turning to R3, note first that
the non-degenerate form $\kappa$ restricts to an $M(F)$-equivariant
perfect pairing $$\kappa:\wlie_{-1} \times \wlie_1 \to F$$ and hence an
$M(F)$-equivariant perfect pairing $$\mu=\bwedge^d \kappa:\bwedge^d
\wlie_{-1} \times \bwedge^d \wlie_1 \to F,$$ where $d = \dim \wlie_1$
Note that when $d=0$, $\mu$ is the multiplication pairing $F
\times F \to F$. Also note for $m \in M(F)$ that $(\bwedge^d \Ad)(m)$
acts on $\bwedge^d \wlie_{\pm 1}$ as multiplication with
$\det(\Ad(m)_{\mid \wlie_1})^{\pm 1}.$ Let $0\not = \tau \in \bwedge^d \wlie_{-1}$.
\begin{prop} (compare \cite{rao-orbital}*{Lemma 2})
\label{prop:phi-prop}
 Define the function $\phi:\wlie_2 \to \R_{\ge 0}$ by the rule
  \begin{equation*}
    \phi(Y) = |\mu(\tau,\bwedge^d(\ad(Y))\tau) |^{1/2}.
  \end{equation*}
  Then $\phi(X) > 0$, and for each $m \in M(F)$ and $Y \in \wlie_2$,
  we have 
  \begin{equation*}
      \phi(\Ad(m)Y) = |\det(\Ad(m)_{\mid \wlie_1}|\phi(Y).
  \end{equation*}
\end{prop}

\begin{proof}
  By Proposition \ref{prop:kr} one knows that $C_G(X) \subset P$; thus
  $\Lie(C_G(X)) \subset \Lie(P)$. Since the orbit of $X$ is separable,
  one knows that $\clie_\glie(X) = \Lie(C_G(X)) \subset \Lie(P)$. This
  implies that $\clie_\glie(X) \cap \glie(-1) = 0$, and so
  $\ad(X):\glie(-1) \to \glie(1)$ is bijective. It follows that
  $\phi(X) >0$.
  
  We have the identity $\ad(\Ad(m)Y) = \Ad(m) \circ \ad(Y) \circ
  \Ad(m^{-1}):\wlie_{-1} \to \wlie_1$.  Thus
  $\bwedge^d(\ad(\Ad(m)Y))\tau = \det(\Ad(m)_{\mid \wlie_1})^2
  \bwedge^d(\ad(Y))\tau$. This implies the second assertion.
\end{proof}
 
This proposition verifies R3, and in view of what was said before,
completes the proof of Theorem \ref{theorem:nilpotent-convergence}.

\subsection{Jordan decomposition}
\label{sec:Jordan-decomposition}

Let $A$ be a linear algebraic group. If $x \in A$ recall that the
Jordan decomposition of $x$ is the expression $x=su$ with $s \in A$
semisimple, $u \in A$ unipotent, and $su=us$. It is a basic fact that
each element has a Jordan decomposition, and that $s$ and $u$ are
uniquely determined.  Similar statements hold for the Jordan
decomposition $X=S+N$ for $X \in \Lie(A)$ (where now $N$ is
nilpotent). In this section, we consider the question of when the
Jordan decomposition of $x \in A(F)$ (and of $X \in \Lie(A)(F)$) is
defined over $F$ in the case when $A$ is a reductive group. Of course, if $x=su
\in A(F)$, we have $s \in A(F)$ if and only if $u \in A(F)$.

\begin{prop}
  \label{prop:jordan-rational}
  Suppose that $p > \rank_{\sems}G +1$.
  \begin{enumerate}
  \item Let $g \in G(F)$, and let $g=su$ be the Jordan decomposition
    of $g$ with $s,u \in G(\overline{F})$.
    Then $s,u \in G(F)$.
  \item Let $X \in \glie(F)$, and let $X=S + N$ be the Jordan
    decomposition of $X$ with $S,N \in \glie(F)$.  Then $S,N \in \glie(F)$.
  \end{enumerate}
\end{prop}

\begin{rem}
  Without our assumption on $p$, the proposition is false. Indeed,
  consider the group $G=GL_{p/F}$, let $f \in F[T]$ be a purely
  inseparable irreducible polynomial of degree $p$, and let $g \in
  \GL_p(F)$ and $X \in \lie{gl}_p(F)$ be any elements having
  characteristic polynomial $f$. Then the semisimple part of each of
  these elements is the scalar matrix $\alpha\cdot I$ where $\alpha$
  is the unique root of $f$ in the algebraically closed extension $k$.
  In particular, this semisimple part is not $F$-rational.
\end{rem}

We begin with a few lemmas.

\begin{lem}
  \label{lem:sep-enough}
  Let $A$ be a linear algebraic group defined over $F$.
  \begin{enumerate}
  \item Let $x \in A(F)$ and let $x=su$ be the Jordan decomposition of
    $x$. Then $u \in A(F)$ if
    and only if $u \in A(F_\sep)$. 
  \item Let $X \in \Lie(A)(F)$ and let $X=S+N$ be the Jordan
    decomposition of $X$. Then $N \in A(F)$ if
    and only if $N \in A(F_\sep)$.
  \end{enumerate}
  Put another way, the Jordan decomposition of an element is defined over $F$ if and
  only if it is defined over $F_\sep$.
\end{lem}

\begin{proof}
  We treat the case $x\in A(F)$; the Lie algebra version is similar.
  Suppose that $u \in A(F_\sep)$ and let $\gamma \in \Gal(F_\sep/F)$.
  To see that $u \in A(F)$, it is enough to see that $u'=\gamma(u)=u$.
  But $x=\gamma(x)=\gamma(su)=s'u'$ (where $u'=\gamma(u)$). Since $s'$
  is semisimple and $u'$ is unipotent, and since evidently
  $s'u'=u's'$, the fact that $u'=u$ follows from the unicity of the
  Jordan decomposition of $x$.
\end{proof}

\begin{lem}
  \label{lem:ss-central-case}
  Let $G$ be a semisimple group over $F$.  Let $x \in G(F)$
  have Jordan decomposition $x=su$, and suppose that $s$ is contained
  in the center of $G$. Then $s,u \in G(F)$.
\end{lem}

\begin{proof}
  In view of the previous lemma, we may as well suppose that $F$ is
  separable closed.  Since the center $Z$ of $G$ is a finite
  diagonalizable subgroup, $Z(\overline{F}) = Z(F)$ (recall we are
  assuming $F$ to be separably closed). Since $s \in Z$, it follows
  that $s \in G(F)$ as desired.
\end{proof}

\begin{lem}
  \label{lem:q-power-centralizer}
  Let $A$ be a linear algebraic group over the algebraically closed
  field $k$ and let $x \in A$ be semisimple. Then $C_A(x) = C_A(x^q)$
  for any $q=p^n$. Similarly, let $X \in \Lie(A)$ be semisimple. Then
  $C_A(X) = C_A(X^{[q]})$, where $X \mapsto
  X^{[p]}$ denotes the $p$-operation on $\Lie(A)$.
\end{lem}

\begin{proof}
  Since $A$ has a faithful matrix representation, it suffices to prove
  the lemma for the group $A=\GL(V)$. Morever, the proof in the Lie
  algebra case is not essentially different, so we discuss only the
  case where $x \in A$. Let $\lambda_1,\dots,\lambda_m$ be the
  distinct eigenvalues of $x$ in $k$. Then
  $\lambda_1^q,\dots,\lambda_m^q$ are the eigenvalues of $x^q$, and
  the lemma follows provided that the $\lambda_i^q$ are all distinct.
  If $\lambda_i^q = \lambda_j^q$, then $\lambda_i/\lambda_j$ is a
  $q$-th root of unity. Since $k$ has characteristic $p$, it follows that
  $\lambda_i = \lambda_j$ so that $i=j$ as desired.
\end{proof}

\begin{proof}[Proof of Proposition \ref{prop:jordan-rational}:]
  We first prove (1).  Since $u$
  is unipotent, $u^q=1$ for some $q=p^n$.  Since $su=us$, we have
  $g^q=s^q \in G(F)$. It follows from Lemma
  \ref{lem:q-power-centralizer} that $C_G^o(s)=C_G^o(s^q)= C_G^o(g^q)$.
  Since $g^q$ is $F$-rational and semisimple, $C=C_G^o(s)$ is a connected,
  reductive $F$-subgroup of $G$.
  
  Let $C_1 = (C,C)$ be the derived group of $C$. Then $C_1$ is a
  semisimple subgroup of $(G,G)$, and by Lemma
  \ref{lem:good-for-subgroups} the prime $p$ is \emph{very good} for
  $C_1$.
  
  Let $\overline{C}=C/Z$ be the corresponding adjoint group, and let
  $\pi:C \to \overline{C}$ be the canonical surjection. Since $p$ is
  very good for $C_1$, Lemma \ref{lem:quasisimple-simple-adjoint}
  implies that the restriction $\phi=\pi_{\mid C_1}$ of $\pi$ to $C_1$ is a
  separable isogeny $\phi:C_1 \to \overline{C}$.
  
  Since $p$ is good for $G$, it follows from
  \cite{springer-steinberg}*{3.15} that $u \in C=C_G^o(g^q)$; since $s$
  is contained in a maximal torus of $G$, we have also $s \in C$ so
  that $g \in C(F)$. Moreover, $s$ is central in $C$. Consider the
  element $v=\pi(g) \in \overline{C}(F)$.  It follows from
  \cite{springer-LAG}*{11.2.14} that the fiber $\phi^{-1}(v) \subset
  C_1$ is defined over $F$.  That fiber must therefore contain a point
  rational over $F_\sep$; thus, there is some $w \in C_1(F_\sep)$ with
  $\phi(w)=v$.  Let $w=s_1u_1$ be the Jordan decomposition of $w$ in
  $C_1$. An application of Lemma \ref{lem:ss-central-case} shows that
  $u_1,s_1 \in C_1(F_\sep)$.

  We now have $w^{-1}g \in C(F_\sep)$. But $\pi(w^{-1}g)=1$ so that
  $w^{-1}g \in Z(F_\sep)$. It follows that $u=u_1 \in C(F_\sep).$ This
  shows that the Jordan decomposition $g=su$ is defined over $F_\sep$.
  It now follows from Lemma \ref{lem:sep-enough} that $s,u \in C(F)$ as
  desired; this proves (1).
  
  The proof of (2) is similar, though a bit easier. Let $X=S+N$ be the
  Jordan decomposition, and again find $n$ large enough so that
  $N^{[q]}=0$ where $q=p^n$. Since $[S,N]=0$, we have
  $X^{[q]}=S^{[q]}$.  Since $X^{[q]} \in \glie(F)$ is semisimple, its
  centralizer $C=C^o_G(X^{[q]})$ is a reductive $F$-subgroup. Again
  let $C_1 = (C,C)$.  Arguing as before, one sees that the
  characteristic is very good for $C_1$.  Since $\Lie(C_1)$ has no
  trivial submodules, one finds that $\Lie(C) = \Lie(C_1) \oplus
  \zlie$ where $\zlie$ is the Lie algebra of the center of $C$. It
  follows that $N \in \Lie(C_1)$ and $S \in \zlie$.  The center of $C$
  is defined over $F$ (e.g. since it is the kernel of the
  $F$-homomorphism $C \to C_{1,\operatorname{adj}}$). Thus, $\zlie$ is
  defined over $F$. Since also $\Lie(C_1)$ is defined over $F$, we
  deduce $\Lie(C)(F) = \Lie(C_1)(F) \oplus \zlie(F)$.  Since $X \in
  \Lie(C)(F)$, it follows that $N \in \Lie(C_1)(F)$ and $S \in \zlie(F)$;
  the proof is now complete.
\end{proof}

\subsection{General orbital integrals on the Lie algebra}
We will now use Ranga Rao's argument \cite{rao-orbital} to deduce the
convergence of a general orbital integral in favorable cases.

\begin{theorem}
  \label{theorem:general-adjoint}
  Let $G$ be an $F$-standard reductive group, let $X \in \glie(F)$
  have Jordan decomposition $X=S + N$.  If $p>\rank_{\sems}G + 1,$
  then $I_X(f)< \infty$ for each $f \in \cc(\glie(F))$.
\end{theorem}

\begin{proof}[Sketch]
  The proof is the same as that of Theorem 2 in \cite{rao-orbital}.
  We outline the argument for the reader's convenience; for full
  details, refer to \emph{loc. cit.}
  
  We note first that Theorem \ref{theorem:nilpotent-convergence}
  remains valid even when the reductive group $G$ is not connected.
  This follows from the fact that
  the $G(F)$ orbit of $X$ is the disjoint union of finitely many
  $G^o(F)$ orbits.
  
  So fix $f \in \cc(\glie(F))$, consider the reductive $F$-group
  \footnote{In fact, under our assumptions on $G$, the centralizer
    $H=C_G(S)$ will also be connected -- see e.g
    \cite{springer-steinberg}*{3.19}. However, we need to apply this
    argument for a proof of Theorem \ref{theorem:general-class} below;
    in that setting the centralizer of the semisimple part of $x$ will
    in general be disconnected, so that the argument described here is
    indeed necessary.}  $H=C_G(S)$ and note that $C=C_G(X)$ is the
  centralizer in $H$ of $N$. Since $N \in \Lie(H)(F)$, it follows from
  Theorem \ref{theorem:nilpotent-convergence}, and the preceeding
  remarks, that we may define for $y \in G(F)$
  \begin{equation}
    \label{eq:function-g}
    g(y) = \int_{H(F)/C(F)} f(\Ad(y)(S + \Ad(h)N)dh^*
  \end{equation}
  where $dh^*$ denotes the invariant measure on $H(F)/C(X)(F)$.
  Then $g$ is continuous in $y$, satisfies $g(yh) = g(y)$ for $h \in
  H(F)$, and the argument in \emph{loc. cit.}
                                %%%%%%%%%%%%%%%%
  \footnote{It is assumed in \cite{rao-orbital} that $G$ is
    semisimple; the argument that the function $g$ has compact support
    given in \emph{loc. cit.} uses the fact that the adjoint
    representation is faithful.  However, it is clear that one can use
    just any faithful linear representation of $G$, rather than the
    adjoint representation.}
                                %%%%%%%%%%%%%%%%
  shows that $g$ has compact support in $G(F)/H(F)$.
  Thus,
  \begin{align*}
    \int_{G(F)/H(F)} g(y)dy^* &= \int_{G(F)/H(F)} dy^* \int_{H(F)/C(F)} 
      f(\Ad(yh)Z)dh^* \\ &= \int_{G(F)/C(F)} f(\Ad(x)Z))dx^* = I_X(f)
  \end{align*}
  is finite.
\end{proof}

Notice that the proof only uses the assumption made on $p$ to know
that $S,N \in \glie(F)$, i.e. that the Jordan decomposition of $X$ is
defined over $F$.

\subsection{Strongly standard groups}

We are going to prove the analogue for groups of Theorem
\ref{theorem:general-adjoint}; to do this, we require a somewhat
stronger hypothesis on our reductive $F$-group $G$. We now explain
this hypothesis. The field $F$ is arbitary.

Consider $F$-groups $H$ which are direct products
\begin{equation*}
  (*)\quad   H = H_1 \times S,
\end{equation*}
where $S$ is an $F$-torus and $H_1$ is a connected, semisimple
$F$-group for which the characteristic is very good. We say that the
reductive $F$-group $G$ is \emph{strongly standard} if there exists a
group $H$ of the form $(*)$ and a separable $F$-isogeny between $G$
and an $F$-Levi subgroup of $H$.  Thus, $G$ is separably isogenous to
$M=C_H(S_1)$ for some $F$-subtorus $S_1 < H$; note that we do not
require $M$ to be the Levi subgroup of an $F$-rational parabolic
subgroup.  It is checked in \cite{mcninch-optimal}*{Proposition 2}
that a strongly standard $F$-group $G$ is $F$-standard in the sense of
\S \ref{sec:general-reductive} of this paper.  Note that any $F$-form of
$\GL_n$ is strongly standard (see Remark 2 of \emph{loc. cit.}) but
that $\SL_n$ is strongly standard just in case $(n,p)=1$.

\subsection{An algebraic analogue of the logarithm}
\label{sub:alg-log}

In characteristic 0, the convergence of unipotent orbital integrals
(on the group) is deduced by Ranga Rao in \cite{rao-orbital} using the
exponential map from the Lie algebra to the group; of course, the
exponential of a nilpotent element is always meaningful in this
setting, and the existence of an open neighborhood of the nilpotent
set on which the exponential converges is also required in \emph{loc.
  cit.} When the characteristic of $F$ is positive, the usual
exponential map may well define an isomorphism between the nilpotent
set and the unipotent set (at least if $p$ is large)
but this isomorphism will never extend to an open neighborhood of the
nilpotent set in $\glie(F)$: the naive exponential of a semisimple
element will never be defined.

To correct this problem, we require a construction used by Bardsley
and Richardson.  For the remainder of \S \ref{sub:alg-log}, $F$ may be
an arbitrary field of characteristic $p$.

\begin{theorem}
  \label{theorem:alg-log}
  Suppose that $H_1$ is a simply connected semisimple $F$-group in
  very good characteristic, and that $G$ is an $F$-Levi subgroup of
  $H_1 \times S$ for some $F$-torus $S$. Let $\UU \subset G$ be the
  unipotent variety, and let $\NN \subset \glie$ be the nilpotent
  variety.  Then there are $G$-stable, $F$-open sets $U \subset G$ and
  $V \subset \glie$ such that $\UU \subset U$ and $\NN \subset V$, and
  there is a $G$-equivariant morphism $\Lambda:G \to \glie$ such that
  \begin{enumerate}
  \item $\Lambda$ is defined over $F$,
  \item $\Lambda_{\mid \UU}:\UU \to \NN$ is an isomorphism of varieties, and
  \item $\Lambda_{\mid U}:U \to V$ is surjective and {\'e}tale.
  \end{enumerate}
\end{theorem}

We will first prove a technical result.
%  If $(\rho,W)$ is a
%\emph{semisimple} representation of an algebraic $F$-group defined
%over $F$, write
%\begin{equation*}
%  W \iso \bigoplus_{i=1}^d W_i
%\end{equation*}
%where the $W_i$ are the \emph{isotypic} components.  Write $T_W \le
%\GL(W)$ for the $d$-dimensional $F$-split torus consisting of all $t
%\in \GL(W)$ for which $t_{\mid W_i}$ is a scalar for each $1 \le i \le
%d$.

\begin{lem}
  Let $H_1$ be a simply connected semisimple $F$-group in very good
  characteristic. Then there is a semisimple $F$-representation
$(\rho,W)$  with the properties:
  \begin{enumerate}
  \item[BR1.] $d\rho:\lie{h}_1 \to \lie{gl}(W)$ is injective, and
  \item[BR2.] there is an $H_1$-invariant $F$-subspace $\mlie \subset
    \lie{gl}(W)$ such that $\lie{gl}(W) = \mlie \oplus
    d\rho(\lie{h}_1)$ and such that $1_W \in \mlie$.
  \end{enumerate}
\end{lem}

\begin{proof}
  If $(\rho,W)$ is a semisimple $F$-representation of $H_1$, BR2 is a
  consequence of BR3:
  \begin{itemize}
  \item[BR3.]  The trace form $\kappa(X,Y) = \tr(d\rho(X) \circ
    d\rho(Y))$ on $\lie{h}_1$ is non-degenerate.
  \end{itemize}
  Indeed, the trace form on $\lie{gl}(W)$ is non-degenerate, and if
  BR3 holds, the first condition of BR2 holds with $\mlie =
  d\rho(\hlie_1)^\perp$. Since $H_1$ is semisimple, $d\rho(\hlie_1)$
  lies in $\lie{sl}(W)$. Thus, $1_W$ is orthogonal to $d\rho(\hlie_1)$
  under the trace form and so lies in $\mlie$.
  
  When $H_1$ is \emph{split}, it follows from
  \cite{springer-steinberg}*{I.5.3} that there is a suitable
  semisimple $F$-representation for which BR1 and BR3 (and hence
  BR2) hold.
  
  In general, we may choose a finite separable extension $F \subset E$
  which splits $H_1$. The preceeding discussion yields an $E$-representation
  $(\rho,W)$ satisfying BR1 and BR3.
  
  By the adjoint property of the restriction of scalars functor, the
  $E$-homomorphism $\rho:H_{1/E} \to \GL(W)$ yields an $F$-homomorphism
  $\rho':H_{1/F} \to R_{E/F}\GL(W)$; the latter group is a closed
  $F$-subgroup of $\GL(W_F)$, where $W_F = R_{E/F}(W)$ denotes the
  $E$-vector space $W$ regarded as an $F$-vectorspace. Thus we may
  regard $\rho'$ as an $F$-representation $(\rho',W_F)$ of $H_1$.
  
  We note that the $F$-representation $(\rho',W_F)$ is semisimple.
  Indeed, extending scalars, there is an isomorphism of
  $H_{1/E}$-representations
  \begin{equation*}
    (\rho' \tensor 1_E,W_F \tensor_F E) \iso \bigoplus_{j=1}^e (\rho,W),
  \end{equation*}
  where $e=[E:F]$; since this scalar extension yields a semisimple
  representation, the original representation was already semisimple.

  If $\phi:W \to W$ is any $E$-linear map, we have
  \begin{equation*}
    \tr_{E/F}(\tr_E(\phi;W)) = \tr_F(\phi;W_F)
  \end{equation*}
  where $\tr_{E/F}:E \to F$ denotes the trace of the separable field
  extension $E/F$.  If $\kappa'$ is the form on $\hlie_1$ determined
  by $\rho'$, this shows that $\kappa' = \tr_{E/F} \circ \kappa$ on
  $\hlie_1(F)$; since $\tr_{E/F}$ is non-0, $\kappa'$ is nondegenerate
  on $\hlie_1(F)$ and hence nondegenerate on $\hlie_1$. This completes the
  proof of the lemma.
\end{proof}

\begin{proof}[Proof of Theorem \ref{theorem:alg-log}]
  The previous lemma gives a semisimple $F$-representation $(\rho,W)$
  of $H_1$ satisfying BR1, BR2 and BR3; we regard $\rho$ as a
  representation of $H$ with $\rho(S)=1$.
        
  We may now define a map $\Lambda:H \to \hlie$ as follows.  For $h
  \in H$, write $\rho(h) = (X,Y) \in \hlie \oplus \mm$ and put
  $\Lambda(h) = X$. Evidently $\Lambda$ is defined over $F$. Since
  $\mm$ is $H$-invariant, $\Lambda$ is $H$-equivariant. Since $1_W \in \mlie$ by BR2,
  $\Lambda(1) =0$.

  The fact that $\Lambda$ satisfies condition (2) of the statement of
  the theorem follows from Corollary 9.3.4 of \cite{BR85}; condition
  (3) follows from Theorem 6.2 in \emph{loc.  cit.}  (``Luna's
  Fundamental Lemma''). This proves the theorem in case $H=G$.
  
  To prove the result for $G$, recall that $G = C_H(S_1)$ for some
  $F$-torus $S_1 \le H$. Thus $\glie = \lie{c}_\hlie(S_1)$ and it is
  clear that that $\Lambda_{\mid G}:G \to \glie$ satisfies conditions
  (1),(2) and (3) of the conclusion of the theorem.
\end{proof}

\begin{rem}
  Note that the group $G$ in the statement of Theorem
  \ref{theorem:alg-log} is strongly standard. It is not clear to the
  author whether the theorem holds more generally for any strongly
  standard group, however.  It holds for instance whenever $G=H$ is a
  semisimple group in very good characteristic such that the trace
  form of the adjoint representation (``Killing form'') is
  non-degenerate. However, this latter condition is not always true;
  for instance, the trace form of the adjoint representation of
  $\operatorname{PSp}(V)$ is identically zero if $p \mid \dim V$.
\end{rem}

\begin{rem}
  The existence of an equivariant $F$-isomorphism $\UU \iso \NN$
  permits us to transfer to $\UU$ a number of the results obtained in
  this paper for nilpotent elements. If $u = \Lambda^{-1}(X)$ for $X
  \in \NN(F)$, then $C_G(u) = C_G(X)$. Moreover, the conjugacy class
  of $u$ is separable if and only if that is so of the orbit of $X$.
  In particular, it follows from Theorem \ref{theorem:split-u-radical}
  that the unipotent radical of $C_G(u)$ is $F$-split under the
  hypothesis that $F$ is perfect or the conjugacy class of $u$ is
  separable.  In case all unipotent classes are separable and $F$ is
  complete for a non-trivial discrete valuation with finite or
  algebraically closed residue field, it follows from Theorem
  \ref{theorem:finiteness} that there are only finitely many
  $G(F)$-orbits on $\UU(F)$.
  
  Note that the Bardsley-Richardson map $\Lambda$ is not necessary; a
  result of T. Springer  allows one to obtain an equivariant $F$-isomorphism
  $\UU \iso \NN$ under milder hypotheses.
\end{rem}

\subsection{Convergence of unipotent orbital integrals}
\label{sub:unipotent}

We now specialize again to the case where $F$ is complete for a
non-trivial discrete valuation and has finite residue field.  Let $G$
be a strongly standard $F$-group.
We are going to prove the following theorem:

\begin{theorem}
  \label{theorem:unipotent-conv}
  Let $u \in G(F)$ be unipotent. Then $I_u(f) < \infty$ for all $f \in
  \cc(G(F))$.
\end{theorem}

We first suppose that $\hat G$ is a second strongly standard $F$-group
and that $\pi: \hat G \to G$ is a separable isogeny.  If $f \in
\cc(G(F))$, we may define $\pi^*(f)$ by the rule:
\begin{equation*}
  \pi^*(f)(g) = f(\pi(g))
\end{equation*}
for $g \in \widehat G(F)$.
Let $\widehat \UU$ and $\widehat \NN$ be the unipotent and nilpotent
varieties for $\widehat G$, and let $\UU$ and $\NN$ be those for $G$.

\begin{lem}
  \label{lem:upper-*}
  \begin{itemize}
  \item[(a)]   Let $f \in \cc(G(F))$. Then $\pi^*(f) \in \cc(\widehat G(F))$.
  \item[(b)] If $u \in G(F)$ there is a unique $\hat u \in \widehat
    G(F)$ such that $u = \pi(\hat u)$. Moreover, let $f \in \cc(G(F)$. Then
    $I_u(f) < \infty$ if and only if $I_{\hat u}(\pi^*(f)) < \infty$.
  \end{itemize}
\end{lem}

\begin{proof}
  (a) The map $\pi^*(f)$ is evidently a continuous function on $\widehat
  G(F)$.  The support of $\pi^*(f)$ is the inverse image under $\pi$
  of the support of $f$; since $\pi$ is an open mapping and since $f$ has
  compact support, this inverse image is contained in a compact set.
  
  (b) The maps $\pi_{\mid \widehat \UU}:\widehat \UU \to \UU$ is an
  equivariant $F$-isomorphism; see e.g.
  \cite{Mc:sub-principal}*{Lemma 27}; it is then clear that in fact
  $I_u(f) = I_{\hat u}(\pi^*(f))$.
\end{proof}

Now suppose that $G$ is an $F$-Levi subgroup of $H=H_1 \times S$ where
$H_1$ is a simply connected semisimple $F$-group in very good
characteristic, and $S$ is an $F$-torus. Write $\UU$ and $\NN$ for the
unipotent and nilpotent varieties for $G$, and denote by $\Lambda:G
\to \glie$ the equivariant $F$-morphism given by Theorem
\ref{theorem:alg-log}. In particular, let $U$ and $V$ be as in the
statement of that theorem.

Since the \'etale map $\Lambda_{\mid U}$ has finite fibers, one may define 
\begin{equation*}
  \Lambda_*(f)(X) =       \sum_{y \in \Lambda^{-1}(X)} f(y)
\end{equation*}
for any function $f \in \cc(G(F))$ whose support is contained in
$U(F)$, and for any $X \in \glie(F)$.

\begin{lem}
  \label{lem:lower-*}
  Let $f \in \cc(G(F))$, and suppose the support of $f$ is
  contained in $U(F)$.  Then $\Lambda_*(f) \in \cc(\glie(F))$.
\end{lem}

\begin{proof}
  The support of $\Lambda_*(f)$ is contained in the image under $\Lambda$ of
  the support of $f$, hence $\Lambda_*(f)$ is compactly supported. The fact
  that $\Lambda_*(f)$ is continuous follows from the inverse function theorem
  \cite{Serre-LieGroups}*{LG3.9}.
\end{proof}

\begin{proof}[Proof of Theorem \ref{theorem:unipotent-conv}]
  It follows from definitions that $G$ is separably isogenous to an
  $F$-group $\widehat G$ which is an $F$-levi subgroup of a group $H =
  H_1 \times S$ where $H_1$ is a semisimple, simply connected
  $F$-group in very good characteristic, and where $S$ is an
  $F$-torus.  According to Lemma \ref{lem:upper-*}, our theorem will
  follow if it is proved for $\widehat G$; thus we replace $G$ by
  $\widehat G$.
  
  We may now find an equivariant $F$-morphism $\Lambda:G \to \glie$ as
  in Theorem \ref{theorem:alg-log}.  With notation as in that theorem,
  one knows that the closure of the class $\Int(G(F))u$ is contained in
  $U(F)$. Thus, it suffices to consider only
  those $f$ whose support is contained in $U(F)$. Let $X = \Lambda(u)
  \in \glie(F)$. Thus $X$ is nilpotent and $\Lambda$ defines an
  isomorphism between $\Int(G(F))u$ and $\Ad(G(F))X$. By Lemma
  \ref{lem:lower-*} and Theorem \ref{theorem:nilpotent-convergence}
  $I_X(\Lambda_*(f)) < \infty$.  But it is clear that $I_u(f) =
  I_X(\Lambda_*(f))$, so the theorem is proved.
\end{proof}

\subsection{Convergence of orbital integrals on $G(F)$}

\begin{theorem}
    \label{theorem:general-class}
    Let $G$ be a strongly standard reductive $F$-group,  and
    assume that $p > \rank_{\sems}G +1$. Let $x \in G(F)$ and let $f
    \in \cc(G(F))$.  Then $I_x(f) < \infty$.
\end{theorem}
\begin{proof}[Sketch]
  This is deduced from Theorem \ref{theorem:unipotent-conv} using
  Ranga Rao's argument \cite{rao-orbital} in the same way that Theorem
  \ref{theorem:general-adjoint} was deduced from Theorem
  \ref{theorem:nilpotent-convergence}. 
\end{proof}

Just as in Theorem \ref{theorem:general-adjoint}, the proof only uses
the assumption made on $p$ to know that the Jordan decomposition of
$x$ is defined over $F$.

\subsection{Convergence for more general groups}
One would hope that the hypothesis $p > \rank_{\sems}G +1$ made in
Theorems \ref{theorem:general-adjoint} and \ref{theorem:general-class}
is unnecessary. It is indeed unnecessary in the following cases. Let
$A$ be a central simple $F$-algebra, and let $G$ be the $F$-group
$\GL_1(A)$, so that $G$ is an inner form of the group
$\GL_{n/F}$ where $\dim_F A = n^2$. The convergence of 
arbitrary orbital integrals is established in \cite{DKV} following
arguments of R. Howe. See also \cite{LaumonI}*{Chapter 4} for a
detailed argument in the case $G=\GL_{n/F}$.

%The technical obstacle to a proof of the convergence of a general
%orbital integral for a standard group $G(F)$ is the failure of the
%Jordan decomposition to be defined over $F$. 

\newcommand{\myresetbiblist}[1]{%
  \settowidth{\labelwidth}{\def\thebib{#1}\BibLabel}%
  \setlength\labelsep{1mm}
  \setlength\leftmargin\labelwidth
  \addtolength\leftmargin\labelsep
}

\newcommand\mylabel[1]{#1\hfil}

%\bibliography{MathBib,books,Preprints}

\end{document}